\input amstex
\input emslab
\input epsf
\magnification=1200
\documentstyle{grgppt01}
\hfuzz=3mm
\parskip=2pt plus 1pt minus 1pt
\footline={\hfil\eightrm \copyright\ Geoffrey Grimmett, \today}

\def\ccw{\overline{\text{\rm co}\,\sW_{p,q}}}
\def\oHH{\overline\HH}
\def\bL{{\boldsymbol\Lambda}}
\def\ZI{Z_{\text{I}}}
\def\ZP{Z_{\text{P}}}

\def\what#1{\widehat{#1}}
\def\SLE{SLE}
\def\rad{\text{\rm rad}}
\def\definition#1#2\par{\subheading{Definition #1}#2\par\smallskip}
\def\g{\gamma}
\def\o{\text{\rm o}}
\def\O{\text{\rm O}}
\def\OP{({\bf OP})}

\def\pd{\partial}
\def\pcg{p_{\text{g}}}
\def\psd#1{p_{\text{\rm sd}}(#1)}

\def\td{{\text{\rm d}}}
\def\pc{p_{\text{\rm c}}}
\def\wpc{\what p_{\text{\rm c}}}
\def\La{\Lambda}
\def\Lpq{\Lambda,p,q}
\def\rc{random-cluster}

\def\bigmid{\,\big\vert\,}
\def\oo{\infty}
\def\ci#1{\citeto{#1}}

\def\es{{\varnothing}}
\def\section#1{\setcounter{Thm}{0}\bigskip\goodbreak\centerline{{\bf\ignorespaces #1}}
  \nobreak\medskip\nobreak\noindent\ignorespaces}
\def\endsection{\medskip\flushpar}
\def\subsection#1#2{\medskip
  \flushpar{\smc#1\quad#2}\nobreak\smallskip\nobreak\flushpar\ignorespaces
  \nobreak\ignorespaces}
\def\ssubsection#1#2{\smallskip\nobreak
  \flushpar{#1.\quad{\it #2\/}}\flushpar\ignorespaces}

\def\sec#1#2{\smallskip\ii{\bf#1.}\hskip7pt \ignorespaces{\bf#2}\par\smallskip}
\def\subs#1#2{\ii{}\ignorespaces\hskip7pt \hbox to 6mm{#1\hfil}\ #2}

\def\capt#1#2{\baselineskip=10pt\botcaption{\baselineskip=10pt\eightpoint{\it Figure #1.}\q
  #2}\endcaption\par}
\def\mletter#1#2#3{\hskip#2cm\lower#3cm\rlap{$#1$}\hskip-#2cm}
\def\lastletter#1#2#3{\hskip#2cm\lower#3cm\rlap{$#1$}\hskip-#2cm\vskip-#3cm}
\def\figure#1\par{\parindent=0pt
  \vbox{\baselineskip=0pt \lineskip=0pt
  \line{\hfil}
  #1}}


\def\lest{\le_{\text{st}}}
\def\lra{\leftrightarrow}
\def\nlra{\nleftrightarrow}
\def\fpq{\phi_{p,q}}

\def\ii{\itemitem}
\def\iii{\par\indent\hangindent3\parindent\textindent}

\def\l{\lambda}
\def\lac{\l_{\text{\rm c}}}
\def\q{\quad}
\def\qq{\qquad}

\def\la{\langle}
\def\ra{\rangle}
\def\up{\uparrow}
\def\down{\downarrow}
\def\s{\sigma}
\def\d{\delta}
\def\b{\beta}
\def\a{\alpha}
\def\t{\theta}
\def\Om{\Omega}
\def\om{\omega}
\def\R{{\Bbb R}}
\def\LL{{\Bbb L}}
\def\PP{{\Bbb P}}
\def\Z{{\Bbb Z}}

\def\EE{{\Bbb E}}
\def\HH{{\Bbb H}}
\def\lra{\leftrightarrow}
\def\sR{{\Cal R}}
\def\sD{{\Cal D}}
\def\sT{{\Cal T}}

\def\sF{{\Cal F}}
\def\sW{{\Cal W}}

\def\tD{{\text{d}}}
\catcode`\@=11
 \def\logo@{}
\catcode`\@=13
\def\today{\number\day
\space\ifcase\month\or
  January\or February\or March\or April\or May\or June\or
  July\or August\or September\or October\or November\or December\fi
  \space\number\year}
\redefine\qed{{\hfill$\square$}}

\topmatter
\title The Random-Cluster Model
\endtitle
\author Geoffrey Grimmett
\endauthor
\address
Statistical Laboratory, University of Cambridge, Wilberforce Road,
Cambridge CB3 0WB, United Kingdom
\endaddress
\email g.r.grimmett{\@}statslab.cam.ac.uk\endemail
\http http://www.statslab.cam.ac.uk/$\sim$grg/ \endhttp


\abstract
The class of random-cluster models is a unification of a 
variety of stochastic processes of significance
for probability and statistical physics, including
percolation, Ising, and Potts models; in addition, 
their study has impact on the theory of certain random
combinatorial structures, and of
electrical networks. Much (but not all) of the physical theory
of Ising/Potts models is best implemented in the context
of the \rc\ representation. 
This systematic summary of \rc\ models 
includes accounts of the fundamental methods and inequalities,
the uniqueness and specification of infinite-volume measures,
the existence and nature of the phase transition, and the 
structure of the subcritical and 
supercritical phases. The theory for two-dimensional lattices
is better developed than for three and more dimensions.
There is a rich collection
of open problems, including some of substantial significance for the general
area of disordered systems, and these are highlighted
when encountered. Amongst the major open questions, there is the problem
of ascertaining the exact nature of the phase transition for general
values of the cluster-weighting factor $q$, and the 
problem of proving that the critical \rc\ model in two dimensions,
with $1\le q\le 4$,
converges when re-scaled to a stochastic L\"owner evolution (\SLE).
Overall the emphasis is upon the \rc\ model
for its own sake, 
rather than upon its applications to Ising and Potts systems.
\endabstract
\endtopmatter

\document

\section{List of contents}
{\parskip=0pt\eightpoint\baselineskip=9.5pt
\sec{1}{Introduction}
\sec{2}{Potts and \rc\ measures}
\subs{2.1}{Random-cluster measures}
\subs{2.2}{Ising and Potts models}
\subs{2.3}{Random-cluster and Ising/Potts coupled}
\subs{2.4}{The limit as $q\downarrow 0$}
\subs{2.5}{Rank-generating functions}
\sec{3}{Infinite-volume \rc\ measures}
\subs{3.1}{Stochastic ordering}
\subs{3.2}{A differential formula}
\subs{3.3}{Conditional probabilities}
\subs{3.4}{Infinite-volume weak limits}
\subs{3.5}{Random-cluster measures on infinite graphs}
\subs{3.6}{The case $q<1$}
\sec{4}{Phase transition, the big picture}
\subs{4.1}{Infinite open clusters}
\subs{4.2}{First- and second-order phase transition}
\sec{5}{General results in $d$ ($\ge 2$) dimensions}
\subs{5.1}{The subcritical phase, $p<\pc(q)$}
\subs{5.2}{The supercritical phase, $p>\pc(q)$}
\subs{5.3}{Near the critical point, $p\simeq\pc(q)$}
\sec{6}{In two dimensions}
\subs{6.1}{Graphical duality}
\subs{6.2}{Value of the critical point}
\subs{6.3}{First-order phase transition}
\subs{6.4}{\SLE\ limit when $q\le 4$}
\sec{7}{On complete graphs and trees}
\subs{7.1}{On complete graphs}
\subs{7.2}{On trees and non-amenable graphs}
\sec{8}{Time-evolutions of \rc\ models}
\subs{8.1}{Reversible dynamics}
\subs{8.2}{Coupling from the past}
\subs{8.3}{Swendsen--Wang dynamics}

}

\section{1. Introduction}
During a classical period, probabilists studied the 
behaviour of {\it independent\/} random 
variables. The emergent theory is rich, and is 
linked through theory and application to 
areas of pure/applied mathematics and to other sciences. 
It is however unable to answer 
important questions from a variety of sources 
concerning large families of {\it 
dependent\/} random variables. Dependence comes 
in many forms, and one of the targets 
of modern probability theory has been to 
derive robust techniques for studying it.
The voice of statistical physics has been 
especially loud in the call for rigour in 
this general area.
In a typical scenario, 
we are provided with an infinity of random 
variables, indexed by the vertices of some 
graph such as the cubic lattice, and which have 
some dependence structure governed by the 
geometry of the graph. Thus mathematicians 
and physicists have had further cause to 
relate probability and geometry. One major outcome of
the synthesis of ideas from physics and probability 
is the theory
of Gibbs states, [\ci{Geo}], which is now established as a significant
branch of probability theory.
	
A classic example of a Gibbs state is the (Lenz--)Ising model 
[\ci{Isi}] for a ferromagnet. 
When formulated on the bounded region $B$ 
of the square lattice $\LL^2$, a random 
variable $\sigma_x$ taking values $-1$ and $+1$ is 
assigned to each vertex $x$ of $B$, and 
the probability of the configuration 
$\sigma$ is proportional to $\exp(-\beta H (\sigma))$, 
where $\beta>0$ and the `energy' $H(\sigma )$ is the 
negative of the sum of $\sigma_x\sigma_y$ over 
all neighbouring pairs $x,y$ of $B$. 
This `starter model' has proved extraordinarily 
successful in generating beautiful and 
relevant mathematics,  and has been useful and provocative in the 
mathematical theory of phase transitions and cooperative phenomena
(see, for example, [\ci{FFS}]).

There are many possible generalisations of the 
Ising model in which the $\sigma_x$ may 
take a general number $q$ of values, rather than $q=2$ only.
One such generalisation, the so-called Potts model [\ci{Pot}], 
has attracted 
especial interest amongst physicists, 
and has displayed a complex and varied structure; 
for example, when $q$ is large, it enjoys a 
discontinuous phase transition,
in contrast to the continuous transition 
believed to take place for small $q$. 
Ising/Potts models are the first 
of three principal ingredients in the story of random-cluster 
models. Note that they are `vertex models' 
in the sense that they involve random variables 
$\sigma_x$ indexed by the vertices $x$ of the underlying graph.

The `(bond) percolation model' was inspired 
by problems of physical type, and emerged from the 
mathematics literature of the 1950s [\ci{BrH}, \ci{dVW}]. 
In this model for a porous medium, each 
edge of a graph is declared `open' (to the passage
of fluid) with probability $p$, and `closed' 
otherwise, different edges having independent states. 
The problem is to determine the typical
large-scale properties of connected components 
of open edges, as the parameter $p$ varies. 
Percolation theory is now a mature part of 
probability, at the core of the study of 
random media and interacting systems, and it is 
the second ingredient in the story of 
random-cluster models. Note that bond percolation is an `edge 
model', in that the random variables are 
indexed by the set of edges of the underlying 
graph. [There is a variant termed `site percolation' in which the
vertices are open/closed at random rather than
the edges.]

The third and final ingredient preceded
the first two, and is the theory of electrical 
networks. Dating back at least to the 1847 paper 
[\ci{Kir}] of Kirchhoff, this sets down a 
method for calculating macroscopic properties 
of an electrical network in terms of its 
local structure. In particular, it explains 
the relevance of counts of certain types of 
spanning trees of the graph. In the modern 
vernacular, an electrical network on a graph 
$G$ may be studied via the properties 
of a `uniformly random spanning tree' on $G$
(see [\ci{BLPS01}]).

These three ingredients seemed fairly 
distinct until Fortuin and Kasteleyn discovered, 
around 1970, [\ci{F72a}, \ci{F72b}, \ci{FK}, \ci{KF}], 
that each features in a certain 
way within a family of probability measures 
of `edge models', parameterised by two 
quantities, $p\in [0,1]$ and $q\in (0,\infty )$. [In 
actuality, electrical networks arise as a 
weak limit of such measures.] These models 
they termed `random-cluster models', and 
they developed the basic theory --- 
correlation inequalities and the like --- 
in a series of papers published thereafter. The true 
power of random-cluster models as a mechanism 
for studying Ising/Potts models has 
emerged only gradually over the intervening thirty years.

We note in passing that the genesis 
of the random-cluster model lay in Kasteleyn's 
observation that each of the three 
ingredients above satisfies certain series/parallel laws: any two 
edges in series (or parallel) may be 
replaced by a single edge in such a way that, if the 
interaction function is adapted accordingly, then the 
distributions of large-scale 
properties remain unchanged.

The family of random-cluster 
measures (that is, probability measures which govern \rc\ models)
is not an extension of the Potts measures. The 
relationship is more sophisticated, and is such 
that {\it correlations\/} for Potts models 
correspond to {\it connections\/} in 
random-cluster models. Thus the {\it correlation 
structure\/} of a Potts model may 
be studied via the {\it stochastic geometry\/} of a 
corresponding random-cluster model. The 
intuition behind this geometrical study comes 
often from percolation, of which the 
random-cluster model is indeed an extension.

It turns out that, in many 
situations involving ferromagnetic
Ising/Potts models, the best way forward is 
via the random-cluster model. As examples of this we mention the existence of 
discontinuous phase transitions [\ci{LMMRS}], 
exact computations in two dimensions 
[\ci{Bax}], the verification of the Wulff 
construction for Ising droplets [\ci{CP}], and the 
Dobrushin theory of interfaces and 
non-translation-invariant measures [\ci{GG01}]. As 
a major exception to the mantra `everything 
worth doing for Ising/Potts is done best via 
random-cluster', we remind the reader of the 
so-called random-current expansion for 
the Ising model, 
wielded with effect in [\ci{A82}, \ci{ABF}, \ci{AF}]
and elsewhere. The random-current
method appears to be Ising-specific, and
has enabled a deep analysis of the 
Ising model unparalleled in more general Potts 
systems. (See Section 5.3.)

The primary target of this review is to 
summarise and promote the theory of random-cluster
models for its own sake. 
In doing so, we encounter many results having direct impact 
on Ising/Potts systems, but we shall not stress such connections. 
Some of the theory has 
been discovered several times by apparently 
independent teams; whilst making a serious 
attempt to list key references, we apologise for 
unwitting omissions of which there will 
certainly be a few. The large number of references to
work of the author is attributable in part to the fact that 
he is acquainted with these contributions.

It is a lesson in humility 
to return to the original Fortuin--Kasteleyn 
papers [\ci{F72a}, \ci{FK}], and especially [\ci{F72b}],
where so much of the basic theory 
was first presented.
These authors may not have followed the slickest of routes, but 
they understood rather well 
the object of their study. Amongst the many papers of
general significance since, we highlight: [\ci{ACCN}],
which brought the topic back to the fore; [\ci{ES}],
where the coupling between Potts and \rc\ models was so beautifully
managed; [\ci{G93}], where the \rc\ model was studied systematically
on infinite grids; [\ci{Hag97}], which links the theory
to several other problems of interest in statistical mechanics;
and [\ci{GHM}], where \rc\ models are placed in the perspective
of stochastic geometry as a tool for studying phase transitions.

This review is restricted mostly to core material for \rc\ models
on the nearest-neighbour cubic lattice in a general number $d$ of dimensions.
Only in passing do we mention such subjects as
extensions to long-range systems
[\ci{ACCN}], mean-field behaviour in high dimensions [\ci{KS}],
and mixing properties [\ci{Al98}]. Neither do we stress the impact that
graphical methods of the \rc\ type
have had on a variety of other disordered systems,
such as the Ashkin--Teller model [\ci{AT}, \ci{Hag97}, \ci{PfV}, \ci{SS}, \ci{WiDo}], 
the Widom--Rowlinson model [\ci{CCK}, \ci{CK96}, \ci{GHag}, \ci{Hag97},
\ci{WiRow}],
or on methods for simulating disordered physical systems 
[\ci{ChMa1}, \ci{ChMa2}, \ci{SW}, \ci{Wol89}]. 

Random-cluster methods may be adapted to systems with
random interactions [\ci{ACCN1}, \ci{G99b}, \ci{Hag97}], and even 
to non-ferromagnetic systems of Edwards--Anderson
spin-glass type [\ci{EdAn}, \ci{New94}, \ci{New97}] 
where, for example, they have been used to prove that,
for a given set $\{J_e\}$ of positive or negative interactions, 
uniqueness of the infinite-volume
Gibbs measure for the ferromagnetic system having interactions
$\{|J_e|\}$ implies uniqueness for the original system.

Amongst earlier papers on \rc\ models, the following include a degree
of review material: [\ci{ACCN}, \ci{BBCK}, 
\ci{GHM}, \ci{G93b}, \ci{Hag97}, \ci{L00}].

Notwithstanding the fairly mature theory 
which has evolved, there remain many open 
problems including some of substantial 
significance for the general area. Many of these 
are marked in the text with the acronym {\bf OP}.    
\endsection

\section{2. Potts and random-cluster processes}
We write $\mu(f)$ for the expectation of a 
random variable $f$ under a probability
measure $\mu$.
\subsection{2.1}{Random-cluster measures}
Let $G=(V,E)$ be a finite graph. An edge $e$ having endvertices $x$ and
$y$ is written as $e = \la x,y\ra$.
A random-cluster measure
on $G$ is a member of a certain class
of probability measures on the set of subsets of the edge set $E$.
We take as state space  the set
$\Omega=\{0,1\}^E$, members of which are vectors $\omega=(\omega(e):e\in E)$.
We speak of the edge $e$ as being {\it open\/} (in $\omega$) if
$\omega(e)=1$, and as being {\it closed\/} if $\omega(e)=0$.
For $\omega\in\Omega$, let $\eta(\omega)=\{e\in E: \omega(e)=1\}$
denote the set of open edges, and let $k(\omega)$ be the number
of connected components (or `open clusters') of the graph
$(V,\eta(\omega))$. 
Note that $k(\omega)$ includes a count of isolated
vertices, that is, of vertices incident to no open edge.
We assign to $\Om$ the $\sigma$-field $\sF$ of all its subsets.

A {\it random-cluster measure\/}  on $G$ has two parameters 
satisfying $0\le p \le 1$ and $q>0$, and is the measure
$\fpq$ on the measurable pair $(\Om,\sF)$ given by
$$
\fpq(\omega) = \frac 1Z \biggl\{\prod_{e\in E} p^{\om (e)} (1-p)^{1-\om
(e)}\biggr\} q^{k(\om )},\qq\om\in\Om,
$$
where the `partition function', or `normalising constant', $Z$ is given by
$$
Z=\sum_{\om\in\Om}\biggl\{\prod_{e\in E} p^{\om (e)} (1-p)^{1-\om
(e)}\biggr\} q^{k(\om )}.
$$
This measure differs from product measure only through
the inclusion of the term $q^{k(\om)}$. Note the difference between
the cases $q\le 1$ and $q\ge 1$: 
the former favours fewer clusters, whereas the latter favours many clusters.
When $q=1$, edges are open/closed independently of
one another. This very special case has been studied in detail
under the titles {\it percolation\/}  and {\it random graphs\/};
see [\ci{Bol}, \ci{G99}, \ci{JLR}]. 
Perhaps the most important values of $q$ are the
integers, since the \rc\ model with
$q\in\{2,3,\dots\}$ corresponds, in a  way sketched in
the next two sections, to the Potts model with $q$ local states.
The bulk of this review is devoted to the theory of \rc\ measures
when $q\ge 1$. The case $q<1$ seems to be harder mathematically
and less important physically. There is some interest in the
limit as $q \downarrow 0$; see Sections 2.4 and 3.6.

We shall sometimes write $\phi_{G,p,q}$ for $\fpq$,
when the choice of graph $G$ is to be stressed.
\comment
Samples from \rc\ measures on the box $[-15,15]^2$ of $\Z^2$ with
`toroidal boundary conditions' are presented in Figure 2.1. They
have been obtained via the method of `coupling from  the past'
(see [\ci{PW96}] and Section 8.2).
\endcomment
Samples from \rc\ measures on $\Z^2$ 
are presented in Figure 2.1. 

\def\jump{1.2cm}
\topinsert
\figure
\centerline{\hfil\epsfxsize=4cm \epsfbox{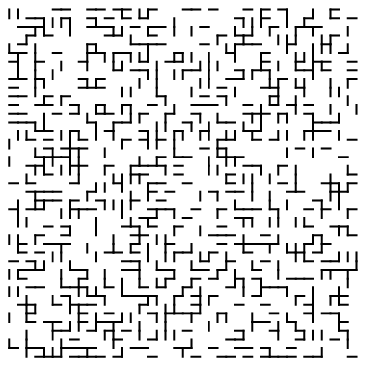}
   \hfil \epsfxsize=4cm \epsfbox{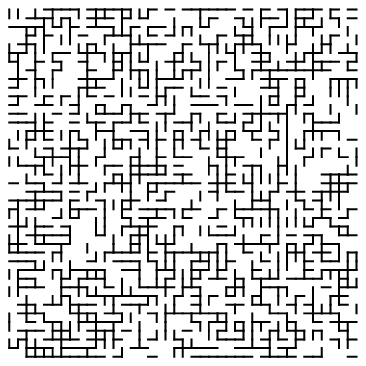}
   \hfil \epsfxsize=4cm \epsfbox{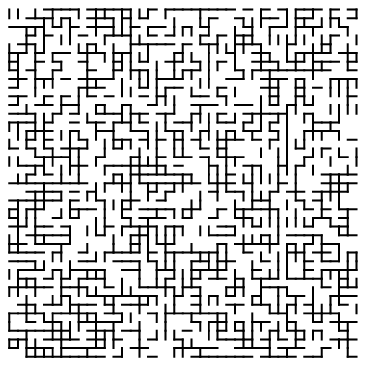}}
\line{\hfil $p=0.30$\hskip\jump \hfil $p=0.45$\hfil\hskip\jump  $p=0.49$\hfil}
\centerline{\hfil\epsfxsize=4cm \epsfbox{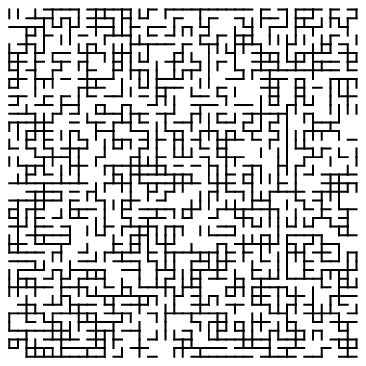}
   \hfil \epsfxsize=4cm \epsfbox{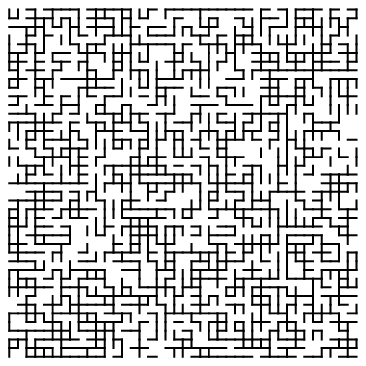}
   \hfil \epsfxsize=4cm \epsfbox{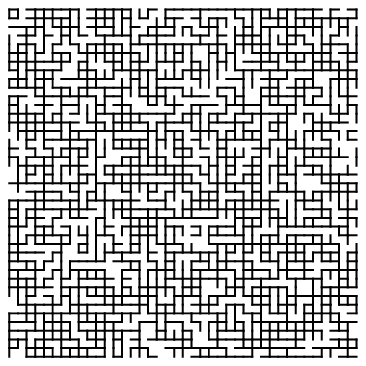}}
\line{\hfil $p=0.51$ \hskip\jump\hfil $p=0.55$ \hfil\hskip\jump $p=0.70$\hfil}

\capt{2.1}{Samples from the \rc\ measure with $q=1$
on the box $[0,40]^2$
of the square lattice. We have set $q=1$ for ease of programming,
the measure being of product form in this case.
The critical value is $\pc(1)=\frac12$.
Samples with more general values of $q$ may be obtained by
the method of `coupling from the past', as described in Section 8.2.}
\endinsert

\comment
\capt{2.1}{Samples from the \rc\ measure with $q=2$
on the box $[-15,15]^2$
of the square lattice, with toroidal boundary conditions, drawn
using the method of `coupling from the past', using programs
kindly provided by
David Wilson.}
\endcomment

\subsection{2.2}{Ising and Potts models}
In a famous experiment, a piece of iron is exposed to a magnetic
field. The field is increased from zero to a maximum, and then
diminishes to zero. If the temperature is sufficiently low, the iron
retains some `residual magnetisation', otherwise it does not. There is
a critical temperature for this phenomenon, often called the {\it Curie
point\/}. The famous (Lenz--)Ising model for such ferromagnetism,  
[\ci{Isi}], may be summarised as follows.
One supposes that particles are positioned at the points of
some lattice in Euclidean space. Each particle may be in either
of two states, representing the physical states of `spin up'
and `spin down'. Spin-values are chosen at random according to
a certain probability measure, known as a {\it Gibbs state\/}, which
is governed by interactions between neighbouring particles.
This measure may be described as follows.

Let $G=(V,E)$ be a finite graph.
We think of each vertex $v\in V$ as being occupied by a particle 
having a random spin.
Since spins are assumed to come in two basic types, we take as sample space
the set $\Sigma =\{ -1,+1\}^V$. The appropriate
probability mass function $\l_{\b,J,h}$ on $\Sigma$ has three
parameters satisfying $0\le \b,J<\infty$ and $h\in\R$, and is given by
$$
\l_{\b,J,h} (\s )=\frac{1}{\ZI}e^{-\b H(\s )},
\qq\s\in\Sigma ,
$$
where the partition function $\ZI$ and the `Hamiltonian' $H :\Sigma\to\R$ 
are given by
$$
\ZI= \sum_{\s\in\Sigma} e^{ -\b H(\s )},\qq
H(\s )=-\sum_{e=\la x,y\ra\in E} J\s_x\s_y-h\sum_{x\in V}\s_x.
$$
The physical interpretation of $\b$ is as the reciprocal $1/T$ of
temperature, of $J$ as the strength of interaction
between neighbours, and of $h$ as the external field.
For reasons of simplicity, we shall consider only the case of
zero external field, and we assume henceforth that $h=0$.
Each edge has equal interaction strength $J$ in the above
formulation. Since $\beta$ and $J$ occur only as a product $\beta J$, 
the measure $\lambda_{\beta,J,0}$ has effectively only a single parameter.
In a more complicated measure not studied here, different edges $e$
are permitted to have different interaction strengths $J_e$.

As pointed out by Baxter, [\ci{Bax}],
the Ising model permits an infinity of generalisations.
Of these, the extension to so-called `Potts models' 
has proved especially fruitful.
Whereas the Ising model permits only two possible spin-values at each
vertex, the Potts model [\ci{Pot}] permits a general number $q\in\{
2,3,\ldots\}$, and is given as follows. 

Let $q$ be an integer satisfying  $q\ge 2$, and take as sample space
$\Sigma =\{ 1,2,\ldots ,q\}^V$. Thus each vertex of $G$ may be
in any of $q$ states. The relevant probability measure is now
given by
$$
\pi_{\b,J,q}(\s )=\frac{1}{\ZP}e^{ -\b H'(\s )},
\qq\s\in\Sigma,
$$
where $\ZP$ is the appropriate normalising constant,
$$
H'(\s )=-\sum_{e=\la x,y\ra}J\d_{\s_x,\s_y},
$$
and $\d_{u,v}$ is the Kronecker delta.
When $q=2$, we have that $\d_{\s_x,\s_y}=\frac12(1+\s_x\s_y)$. It is
now easy to see in this case
that the ensuing Potts model is simply the Ising
model with an adjusted value of $J$.

\subsection{2.3}{Random-cluster and Ising--Potts coupled}
It was Fortuin and Kasteleyn [\ci{F72a}, \ci{F72b}, \ci{FK}, \ci{KF}]
who discovered that Potts models may be
recast as random-cluster models, and furthermore that the relationship
between the two systems facilitates an extended study of
phase transitions in Potts models. Their methods were essentially
combinatorial. In the more modern game, we construct
the two systems on a common probability space, and then 
observe their relationship through their realisations. There
may in principle be many ways to do this, but the standard coupling 
reported in [\ci{ES}] is of special value.

Let $q\in\{2,3,\dots\}$, $0\le p\le 1$,
and let $G=(V,E)$ be a finite graph, as before.  
We consider the product sample space $\Sigma\times\Om$ where
$\Sigma=\{1,2,\dots,q\}^V$ and $\Om=\{0,1\}^E$ as above.
We now define a probability mass
function $\mu$ on $\Sigma\times\Om$ by
$$
\mu (\s ,\om )\propto\prod_{e\in E}\Bigl\{ (1-p)\d_{\om
(e),0}+p\d_{\om (e),1}\d_e(\s )\Bigr\},\qq (\s,\om)\in\Sigma\times\Om,
$$
where 
$\d_e(\s )=\d_{\s_x,\s_y}$ for $e=\la x,y\ra\in E$.
Elementary calculations reveal the following facts.
\ii{(a)} {\it Marginal on $\Sigma$\/}. The marginal measure
$\mu_1 (\s)=\sum_{\om\in\Om}\mu (\s ,\om )$
is given by
$$
\mu_1(\s)\propto\exp\biggl\{\b \sum_e J\d_e(\s )\biggr\}
$$
where $p=1-e^{-\b J}$. This is the Potts measure.

\ii{(b)} {\it Marginal on $\Om$\/}. The second marginal of $\mu$ is
$$
\mu_2(\om )=\sum_{\s\in\Sigma} \mu (\s ,\om )
\propto\biggl\{\prod_e p^{\om (e)}
(1-p)^{1-\om (e)}\biggr\} q^{k(\om )}.
$$
This is the \rc\ measure.

\ii{(c)} {\it The conditional measures\/}. Given $\om$, the
conditional measure on $\Sigma$ is obtained by putting (uniformly)
random spins on entire clusters of $\om$ (of which there are $k(\om
)$). These spins are constant on given clusters, and are independent between
clusters. Given $\s$, the conditional measure on $\Om$ is obtained by
setting $\om (e)=0$ if $\d_e(\s )=0$, and otherwise $\om (e)=1$ with
probability $p$ (independently of other edges).

In conclusion, the measure $\mu$ is a coupling of a Potts measure
$\pi_{\b ,J,q}$ on $V$, together with the random-cluster measure
$\fpq$ on $\Omega$.
The parameters of these measures are related by the equation
$p=1-e^{-\b J}$. Since $0\leq p\leq 1$, this is possible only
if $\b J\geq 0$.

This special coupling may be used in a particularly simple way 
to show that correlations
in Potts models correspond to open connections in random-cluster models.
When extended to infinite graphs, this implies as discussed
in Section 4 that the phase
transition of a Potts model corresponds to the creation of
an infinite open cluster in the random-cluster model. Thus
arguments of stochastic geometry, and particularly those developed for
the percolation model,  may be harnessed directly in order to
understand the correlation structure of the Potts system.
The basic step is as follows.

We write $\{x\lra y\}$ for the set of all $\omega\in\Omega$ for
which there exists an open path joining vertex $x$ to vertex $y$. The complement
of the event $\{x\lra y\}$ is denoted $\{x\nlra y\}$.

The `two-point correlation function' of the
Potts measure $\pi_{\b,J,q}$ on the finite graph 
$G=(V,E)$ is defined to be the function
$\tau_{\b,J,q}$ given by
$$
\tau_{\b,J,q}(x,y)=\pi_{\b,J,q}(\s_x=\s_y)-\frac{1}{q},\qq x,y\in V.
$$
The term $q^{-1}$ is the
probability that two independent and uniformly distributed spins are equal.
The `two-point connectivity function' of the random-cluster measure
$\phi_{p,q}$ is defined as the function $\fpq (x\lra y)$ for $x,y\in V$, 
that is, the probability that $x$ and $y$
are joined by a path of open edges.
It turns out that these `two-point functions'
are (except for a constant factor) the same.

\proclaim{Theorem 2.\Thm{} (Correlation/connection) [\ci{KF}]} 
If $q\in\{ 2,3,\ldots\}$ and $p=1-e^{-\b J}$
satisfies $0\leq p\leq 1$, then
$$
\tau_{\b,J,q}(x,y)=(1-q^{-1})\fpq (x\lra y),\qq x,y\in V.
$$
\endproclaim
\label{corrconn}

\demo{Proof} The indicator function of an event $A$ is denoted $1_A$.
We have that
$$\align
\tau_{\b,J,q}(x,y)&=\sum_{\s ,\om}\Bigl\{ 1_{\{\s_x=\s_y\}} (\s
)-q^{-1}\Bigr\}\mu (\s ,\om )\\
&=\sum_\om\fpq (\om )\sum_\s\mu (\s\mid\om )\Bigl\{
1_{\{\s_x=\s_y\}} (\s )-q^{-1}\Bigr\}\\
&=\sum_\om\fpq (\om )\Bigl\{ (1-q^{-1}) 1_{\{ x\lra y\}} (\om )
+0\cdot 1_{\{ x\nlra y\}} (\om )\Bigr\} \\
&=(1-q^{-1})\fpq (x\lra y),
\endalign
$$
where $\mu$ is the above coupling of the Potts and random-cluster measures.
\qed\enddemo

The theorem may be generalised as follows.
Suppose we are studying the Potts model, and are interested in some
`observable' $f:\Sigma\to\R$.
The mean value of $f(\s )$ satisfies
$$
\align
\pi_{\b,J,q}(f)&=\sum_\s f(\s )\pi_{\b,J,q} (\s )=
   \sum_{\s ,\om} f(\s )\mu (\s ,\om
)\\
&=\sum_\om F(\om )\phi_{p,q} (\om ) =\phi_{p,q}(F)
\endalign
$$
where $F:\Om\to\R$ is given by
$$
F(\om ) =\mu(f\mid \omega)=\sum_\s f(\s )\mu (\s\mid\om ).
$$
The above theorem is obtained in the case
$f(\s )=\d_{\s_x,\s_y}-q^{-1}$, where $x,y\in V$.

The Potts models considered above have zero external field. Some 
complications arise when an external field is added; see the discussions
in [\ci{Al01}, \ci{BBCK}].

\subsection{2.4}{The limit as $q\downarrow 0$}
Let  $G=(V,E)$ be a finite  connected graph,
and let $\fpq$ be the \rc\ measure on the associated sample
space $\Omega=\{0,1\}^E$.
We consider first the weak limit of $\fpq$ as $q\downarrow 0$
for fixed $p\in(0,1)$. 
This limit may be ascertained by observing that the dominant terms
in the partition function
$$
Z(p,q)=\sum_{\omega\in\Omega} p^{|\eta(\om)|}(1-p)^{|E\setminus\eta(\om)|}
q^{k(\om)}
$$
are those for which $k(\omega)$ is a minimum, that is, those with
$k(\omega)=1$. It follows that $\lim_{q\downarrow 0}\fpq$
is precisely the product measure $\phi_{p,1}$ (that is, percolation
with intensity $p$)
conditioned on the
resulting graph $(V,\eta(\om))$ being connected. A more interesting
limit arises if we allow $p$ to converge to 0 with $q$, as follows.

The \rc\ model originated in a systematic
study by Fortuin and Kasteleyn of systems of 
a certain type which satisfy certain
parallel and series laws. Electrical 
networks are the best known such systems ---
two parallel (respectively, series) connections of resistances $r_1$ and
$r_2$ may be replaced by a single connection 
with resistance $(r_1^{-1}+r_2^{-1})^{-1}$
(respectively, $r_1+r_2$). Fortuin and Kasteleyn [\ci{FK}] realised
that the electrical-network theory of a graph $G$ is
related to the limit as $q\downarrow 0$ of the \rc\ model
on $G$. Their argument may be expanded as follows.

Suppose $p=p_q$ is related to $q$ in such a way that
$p\to 0$ and $q/p\to 0$ as $q\to 0$.
We may write $Z(p,q)$ as
$$
Z(p,q)=(1-p)^{|E|} \sum_{\om\in\Om}\left(\frac p{1-p}\right)^{|\eta(\om)|+k(\om)}
\left(\frac{q(1-p)}{p}\right)^{k(\om)}.
$$
Note that $p/(1-p)\to 0$ and $q(1-p)/p\to 0$ as $q\to 0$.
Now $k(\om)\ge 1$ and $|\eta(\om)|+k(\om)\ge |V|$ for all
$\om\in\Om$; these two inequalities are satisfied 
simultaneously with equality 
if and only if $\eta(\om)$ is a spanning tree of $G$.
It follows that, in the limit as $q\to 0$, the `mass' is concentrated
on such configurations, and it is easily
seen that the limit mass is uniformly distributed. 
That is, $\lim_{q\downarrow 0}\fpq$ is a probability measure 
which selects, uniformly
at random, a spanning tree of $G$; in other words, the limit
measure is $\phi_{\frac12,1}$ conditioned on the resulting
graph being a spanning tree.

The link to the theory of electrical networks is now provided by Kirchhoff's
theorem [\ci{Kir}], which expresses effective
resistances in terms of counts of spanning trees. See
also [\ci{Hag95}].

The theory of random spanning trees is beautiful in its own right
(see [\ci{BLPS01}]), and is
linked in an important way to the emerging field of stochastic
growth processes of `stochastic L\"owner evolution' (\SLE) type
(see [\ci{LSW7}, \ci{RS01}]), to which we return in Section 6.4.
Another limit emerges if $p=q$ and $q\downarrow 0$, namely
uniform measure on the set of forests of $G$. 
More generally, take $p=\a q$ where $\a\in(0,\infty)$ is constant,
and take the limit as $q\downarrow 0$. The limit measure is
the percolation measure
$\phi_{\beta,1}$ conditioned on the non-existence of open circuits,
where $\beta=\a/(1+\a)$. If $p/q\to 0$ as $p,q\to 0$,
the limit measure is concentrated on the empty set of edges.

\subsection{2.5}{Rank-generating functions}
The partition functions of Potts and \rc\ measures are
particular evaluations of rank-generating functions, defined as follows.
The {\it rank-generating
function\/} of the simple graph $G=(V,E)$ is the function
$$
W_G(u,v)=\sum_{E'\subseteq E} u^{r(G')} v^{c(G')},\qq u,v\in\R,
$$
where $r(G')=|V|-k(G')$ is the {\it rank\/} of the graph $G'=(V,E')$,
and $c(G')=|E'|-|V|+k(G')$ is its {\it co-rank\/}; here, $k(G')$
denotes the number of components of the graph $G'$. The
rank-generating function has various useful properties, and occurs in
several contexts in graph theory; see [\ci{Big}, \ci{Tut}]. It
crops up in other forms also. For
example, the function
$$
T_G(u,v)= (u-1)^{|V|-1} W_G\bigl( (u-1)^{-1}, v-1\bigr)
$$
is known as the {\it dichromatic\/} (or {\it Tutte\/}) {\it
polynomial\/}, [\ci{Tut}]. The partition function $Z=Z_G$ of the random-cluster
measure on $G$ with parameters $p$, $q$
is easily seen to satisfy
$$
Z_G= q^{|V|} (1-p)^{|E|} W_G
\left(\frac p{q(1-p)}, \frac p{1-p}\right) ,
$$
a relationship which provides a link with other classical quantities
associated with a graph. See [\ci{Big}, \ci{Big77}, 
\ci{F72a}, \ci{WeMe}] also.
\endsection

\section{3. Infinite-volume random-cluster measures}
It is in the infinite-volume limit that \rc\ measures exhibit phase
transitions.  There are two ways of constructing \rc\ measures 
on infinite graphs, namely
by taking weak limits as a finite domain approaches
the infinite system, and by studying measures on
the infinite graph having the
`correct' conditional versions. Such matters are discussed
in this section, which begins with a summary of certain valuable properties
of \rc\ measures on finite graphs.

\subsection{3.1}{Stochastic ordering}
The stochastic ordering of measures provides a technique fundamental
to the study of random-cluster measures. Let $G=(V,E)$ be
a finite or countably infinite graph as above; let $\Om=\{0,1\}^E$,
and let $\sF$ be the $\sigma$-field of $\Omega$ 
generated by the finite-dimensional
cylinders.
Note first that $\Om$ is a partially ordered set with partial order
$\om_1\le\om_2$ if $\om_1(e)\le\om_2(e)$ for all $e$. A random
variable $f:\Om\to\R$ is called {\it increasing\/} if
$f(\om_1)\le f(\om_2)$ whenever $\om_1\le \om_2$.
An event $A\in\sF$ is called {\it increasing\/} if its indicator
function $1_A$ is increasing. The word `decreasing' should be 
interpreted in the natural way. Given two probability measures
$\mu_1$, $\mu_2$ on $\Om$, we write $\mu_1\lest \mu_2$, and
say that $\mu_1$ is stochastically smaller than $\mu_2$,
if $\mu_1(f)\le \mu_2(f)$ for all 
bounded increasing random variables $f$ on $\Om$.

We return now to the case when $G$ is a {\it finite\/} graph.
Let $\mu_1$, $\mu_2$ be probability measures on $\Omega$,
and assume for the moment that the $\mu_i$ are strictly positive 
in the sense that
$\mu_i(\om)>0$ for all $\om\in\Om$. 
An important sufficient condition  for
the inequality  $\mu_1\lest \mu_2$
was found by Holley [\ci{Hol}], namely that
$$
\mu_1(\om_1\vee \om_2)\mu_2(\om_1\wedge\om_2)\ge \mu_1(\om_1)\mu_2(\om_2)
\qq\text{for all } \om_1,\om_2\in\Om,
$$
where $\om_1\vee\om_2$ and $\om_1\wedge\om_2$ are the maximum and minimum
configurations given respectively as $\max\{\om_1(e),\om_2(e)\}$
and $\min\{\om_1(e),\om_2(e)\}$, for $e\in E$.
A probability measure $\mu$ on $\Om$ is said to have
the {\it FKG lattice property\/} if
$$
\mu(\om_1\vee \om_2)\mu(\om_1\wedge\om_2)\ge \mu(\om_1)\mu(\om_2)
\qq\text{for all } \om_1,\om_2\in\Om,
$$
amd it is a consequence of Holley's argument that
any strictly positive measure with the FKG lattice property satisfies
the so-called FKG inequality.  This amounts to the following
for \rc\ measures.

\proclaim{Theorem 3.\Thm\hbox{\ }(FKG inequality) [\ci{F72b}, \ci{FKG}]} 
Suppose that $0\le p\le 1$ and $q\geq 1$. If $f$
and $g$ are increasing functions on $\Om$, then
$\fpq (fg)\geq \fpq(f)\fpq (g)$.
\endproclaim
\label{fkg}

Specialising to indicator functions, we obtain that
$$
\phi_{p,q} (A\cap B)\geq\phi_{p,q} (A)\phi_{p,q} (B)\qq\hbox{for
increasing events $A,B$},
$$
whenever $q\geq 1$. It is not difficult to see that the FKG inequality
does not generally hold when $0<q<1$.

Holley's theorem leads easily to the following comparison inequalities,
which were first proved by Fortuin.

\proclaim{Theorem 3.{\Thm{}}\hbox{\ }(Comparison inequalities) [\ci{F72b}]} 
It is the case that
$$
\align
\phi_{p',q'}&\leq\phi_{p,q}\qq\text{if}\q q'\geq q,\ q'\geq 1,
\text{ and } p'\leq p,\\
\phi_{p',q'}&\geq\phi_{p,q}\qq\text{if}\q q'\geq q,\ q'\geq 1,
\text{ and } \frac{p'}{q'(1-p')} \geq \frac p{q(1-p)}.
\endalign
$$
\endproclaim
\label{dom}

\subsection{3.2}{A differential formula}
One way of estimating the probability of an event $A$
is via an estimate of its derivative $d\fpq(A)/dp$.
When $q=1$, there is a formula for this derivative
which has proved very useful, and which is commonly
attributed to Russo, see
[\ci{BP65}, \ci{G99}, \ci{Ru81}]. This formula may be generalised
to \rc\ measures as follows. The proof is an exercise in the differentiation
of summations.

\proclaim{Theorem 3.\Thm{} [\ci{BGK}]}  Let $0<p<1$, $q>0$, and
let $\fpq$ be the corresponding random-cluster measure on a finite
graph $G=(V,E)$. Then
$$
\frac{d}{dp}\fpq(A)=\frac{1}{p(1-p)}
\bigl\{\fpq(|\eta|1_A)-\fpq(|\eta|)\fpq(A)\bigr\}
$$
for any event $A$, where $|\eta|=|\eta(\om)|=\sum_{e\in E}\om(e)$ 
is the number of open edges of the
configuration $\om$.
\endproclaim
\label{diff}

\subsection{3.3}{Conditional probabilities}
Whether or not an edge $e$ is open
depends on
the configuration on $E\setminus\{e\}$, and a further
important property of \rc\ measures summarises the nature of this dependence.

For $e \in E$, we denote by $G\setminus e$
(respectively, $G.e$) the graph obtained from $G$ by deleting (respectively,
contracting) $e$.  We write $\Om_e = \{0,1\}^{E\setminus\{e\}}$;
for $\om \in \Om$ we define $\om_e \in \Om_e$ by
$\om_e(f) = \om(f)$ for $f \ne e$. 
For $e=\la x,y\ra$, we write $K_e$ for the event that $x$ and $y$
are joined by an open path not using $e$.

\proclaim{Theorem 3.\Thm{} [\ci{F72b}]} 
Let $e\in E$. We have that
$$
\phi_{G,p,q}\bigl(\om\bigmid \om(e)=j\bigr)
=\cases  \phi_{G\setminus e,p,q}
(\om_e)&\hbox{if } j=0,  \\
\phi_{G.e,p,q}
(\om_e) &\hbox{if } j=1, 
\endcases
$$
and
$$
\phi_{G,p,q}\bigl(\om(e)=1\bigmid \om_e\bigr) = \cases p &\text{if } \om_e\in K_e,\\
\dfrac p{p+(1-p)q} &\text{if } \om_e\notin K_e.
\endcases
$$
\endproclaim
\label{tower}

That is to say, the effect of conditioning on the
absence or presence of an edge $e$ is to replace the measure
$\phi_{G,p,q}$ by the \rc\ measure on the respective graph
$G\setminus e$ or $G.e$. Secondly, the conditional
probability that $e$ is open, given the configuration elsewhere, depends only
on whether or not $K_e$ occurs, and is then given by the stated formula.
The proof is elementary. The final equation of the theorem leads
to properties of \rc\ measures referred to elsewhere
as `insertion tolerance' and the `finite-energy property'.

\subsection{3.4}{Infinite-volume weak limits}
In studying random-cluster measures on infinite graphs, we restrict ourselves
to the case of the hypercubic lattice in $d$ dimensions, where $d\geq
2$; similar observations are valid in greater generality. Let $d\geq
2$, and let $\Z^d$ be the set of all $d$-vectors of integers; for
$x\in\Z^d$, we normally write $x=(x_1,x_2,\ldots ,x_d)$. For
$x,y\in\Z^d$, let
$$
\| x-y\| =\sum_{i=1}^d |x_i-y_i|.
$$
We place an edge $\la x,y\ra$ between $x$ and $y$ if and only if $\|
x-y\| =1$; the set of such edges is denoted by $\EE^d$, and we write
$\LL^d=(\Z^d,\EE^d)$ for the ensuing lattice. For any subset $S$ of
$\Z^d$, we write $\partial S$ for its boundary, that is, 
$$
\partial S=\{ s\in S:\la s,t\ra\in\EE^d\ \hbox{for some $t\not\in S\}$}.
$$

Let $\Om=\{0,1\}^{\EE^d}$, and let $\sF$ be the $\s$-field
of subsets of $\Om$ generated by the finite-dimensional cylinders.
The letter $\Lambda$ is used to denote
a finite box of $\Z^d$, which is to say that
$\Lambda =\prod_{i=1}^d [x_i,y_i]$
for some $x,y\in\Z^d$; we interpret $[x_i,y_i]$ as the set $\{
x_i,x_i+1,x_i+2,\ldots ,y_i\}$. The set $\Lambda$ generates a subgraph of
$\LL^d$ having vertex set $\Lambda$ and edge set $\EE_\Lambda$
containing all $\la x,y\ra$ with $x,y\in\Lambda$.

We are interested in the `thermodynamic limit' (as
$\Lambda\up\Z^d)$ of the random-cluster measure on the finite box
$\Lambda$. In order to describe such weak limits, we shall need 
to introduce the notion
of a `boundary condition'.

For $\xi\in\Om$,
we write $\Om_\Lambda^\xi$ for
the (finite) subset of $\Om$ containing all configurations $\om$ satisfying
$\om (e)=\xi (e)$ for $e\in \EE^d\setminus\EE_\Lambda$;
these are the configurations which `agree with $\xi$ off $\La$'. 
For $\xi\in\Om$ and values
of $p$, $q$ satisfying $0\leq p\leq 1$, $q>0$, we define $\phi_{\Lpq}^\xi$
to be the random-cluster measure on the finite graph
$(\Lambda ,\EE_\Lambda )$ `with boundary condition $\xi$'; this is the
equivalent of a `specification' for Gibbs states.
More precisely, let $\phi_{\Lpq}^\xi$
be the probability measure on the pair $(\Om ,\sF )$ given by
$$
\phi_{\Lpq}^\xi (\om ) = \cases\dfrac 1{Z_{\Lpq}^\xi}
\bigg\{\displaystyle\prod_{e\in\EE_\Lambda} p^{\om (e)} (1-p)^{1-\om (e)}\bigg\}
q^{k(\om ,\Lambda )} &\text{if $\om\in\Om_\Lambda^\xi$},\\
0 &\text{otherwise},
\endcases
$$
where $k(\om ,\Lambda )$ is the number of components of the graph
$(\Z^d,\eta (\om ))$ which intersect $\Lambda$, and where
$Z_{\Lpq}^\xi$ is the appropriate normalising constant
$$
Z_{\Lpq}^\xi =\sum_{\om\in\Om_\Lambda^\xi}\bigg\{\prod_{e\in\EE_\Lambda}
p^{\om (e)} (1-p)^{1-\om (e)}\bigg\} q^{k(\om ,\Lambda )}.
$$
Note that $\phi_{\Lpq}^\xi (\Om _\La^\xi) = 1$.

\comment
One works often on a finite box $\La$ of $\Z^d$, and then
takes the limit as $\La\uparrow\Z^d$. This limit is to be interpreted along a
sequence $(\La_n: n\ge 1)$ of boxes such that,
for all $m$, $\La_n\supseteq[-m,m]^d$ for
all large $n$.

\definition{3.\Thm} Let $0\le p\le 1$ and $q>0$.
A probability measure $\phi$ on $(\Om ,\sF
)$ is called a {\it limit random-cluster measure\/} with parameters $p$
and $q$ if there exists $\xi\in\Om$ and an increasing sequence
$(\Lambda_n:n\geq 1)$ of boxes, satisfying $\Lambda_n\to\Z^d$ as
$n\to\infty$, such that
$$
\phi_{\Lambda_n,p,q}^\xi\Rightarrow\phi\q\text{as $n\to\infty$}
$$
where `$\Rightarrow$' denotes weak convergence. The set of all such
measures is denoted by $\sW_{p,q}$.
\label{wldef}
\endcomment

\definition{3.\Thm} Let $0\le p\le 1$ and $q>0$.
A probability measure $\phi$ on $(\Om ,\sF
)$ is called a {\it limit random-cluster measure\/} with parameters $p$
and $q$ if there exist $\xi\in\Om$ and 
a sequence $\bL=(\Lambda_n:n\ge 1)$ of boxes
satisfying $\La_n\to\Z^d$ as $n\to\oo$ such that
$$
\phi_{\La_n,p,q}^\xi \Rightarrow \phi \qq\text{as } n\to \oo.
$$
The set of all such
measures $\phi$ is denoted by $\sW_{p,q}$, and the closed convex
hull of $\sW_{p,q}$ is denoted $\ccw$.
\label{wldef}

In writing $\Lambda_n\to\Z^d$ we mean that,
for all $m$, $\La_n\supseteq[-m,m]^d$ for
all large $n$. The arrow `$\Rightarrow$' denotes
weak convergence.

It might seem reasonable to define a limit
\rc\ measure to be any weak limit of the form
$\lim_{n\to\oo} \phi_{\La_n,p,q}^{\xi_n}$ for some sequence
$(\xi_n:n\ge 1)$ of members of $\Om$ 
and some sequence $\bL=(\Lambda_n:n\ge 1)$ of boxes
satisfying $\La_n\to\Z^d$.
It may however be shown that this adds no extra generality
to the class as defined above, [\ci{G93}].
The dependence
of the limit measure $\phi$ on the choice of sequence $(\Lambda_n)$
can be subtle, especially when $q<1$ \OP.

It is standard that 
$\sW_{p,q}\ne\es$ for all $0\leq p\leq 1$, $q>0$, and 
one way of seeing this is
as follows. The sample space $\Om$ is the product of discrete spaces,
and is therefore compact. It follows that any class of
probability measures on $\Om$ is tight, and hence relatively
compact 
(see the account of Prohorov's theorem in [\ci{Bil}]),
which is to say that any infinite 
sequence of probability measures contains a weakly
convergent subsequence.

When does the limit $\lim_{n\to\infty} \phi_{\Lambda_n,p,q}^\xi$
exist, and when does it depend on the choice of boundary condition $\xi$?  
The
FKG inequality provides a route to a partial answer to this important question.
Suppose for the moment that $q\ge 1$.
Two extremal boundary conditions of special importance are provided
by the configurations $0$ and $1$, comprising 
`all edges closed' and `all edges open'
respectively.
One speaks of configurations in
$\Om_\Lambda^0$ as having `free' boundary conditions, and
configurations in $\Om_\Lambda^1$ as having `wired' boundary
conditions. 

\proclaim{Theorem 3.\Thm{}\hbox{\ }(Thermodynamic limit) 
[\ci{ACCN}, \ci{BC96}, \ci{F72b}, \ci{G93b}, \ci{G93}]} 
Suppose $0\le p\le 1$ and $q\geq 1$.
\flushpar
{\rm (a)} Let
$\bL=(\Lambda_n:n \ge 1)$ be a sequence of boxes satisfying
$\Lambda_n\to\Z^d$ as $n \to \infty$. 
The weak limits
$$
\phi_{p,q}^b =\lim_{n\to\infty} \phi_{\Lambda_n,p,q}^b,
\qq\hbox{for $b=0,1$},
$$
exist and are independent of the choice of $\bL$.
\flushpar
{\rm (b)} We have that each  $\phi_{p,q}^b$ 
is translation-invariant, and
$$
\fpq^0 \lest \phi \lest \fpq^1\qq\text{for all } \phi\in\sW_{p,q}.
$$
\flushpar
{\rm (c)}  For $b=0,1$, the measure $\phi_{p,q}^b$ is
ergodic, in that 
any translation-invariant random variable is $\phi_{p,q}^b$-a.s.\ constant.
\endproclaim
\label{thermo}

The FKG inequality underlies all parts of Theorem
3.\refer{thermo}.
The claim (c) of ergodicity has until recently been considered
slightly subtle (see the discussion after the forthcoming Theorem
3.\refer{rcexist}) but an easy proof may be found in [\ci{L00}].

It follows from the inequality of part (b) that
$|\sW_{p,q}|=1$
if and only if $\phi_{p,q}^0=\phi_{p,q}^1$. It is an important open
problem to determine for which $p$, $q$ this holds, and we
shall return to this question in Section
5 \OP. For the moment, we note one 
sufficient condition for uniqueness,
proved using a certain convexity property of the 
logarithm of a partition function $Z$.

\proclaim{Theorem 3.\Thm{}\ [\ci{G94a}, \ci{G93}]}
Let $q \ge 1$. There exists a subset $\sD_q$ of $[0,1]$,
at most countably infinite in size,
such that $\phi_{p,q}^0=\phi_{p,q}^1$, and hence $|\sW_{p,q}|=1$,
if $p \notin \sD_q$.
\endproclaim
\label{countable}

It is believed but not proved \OP\ that: for any given $q\ge 1$,
$\sD_q$ either is empty
or consists of a singleton (the
critical point,
to be defined in Section 4),
the former occurring
if and only if 
$q$ is sufficiently small.

\subsection{3.5}{Random-cluster measures on infinite graphs}
One may define a class of measures on the infinite lattice without
having recourse to weak limits. 
The following definition of a random-cluster measure
is based upon the Dobrushin--Lanford--Ruelle (DLR)
definition of a Gibbs state, [\ci{Dob}, \ci{Geo}, \ci{LR}].
It was introduced in [\ci{G93b}, \ci{G94a}], and discussed further in 
[\ci{BC96}, \ci{G93}].
For any box $\La$, we write $\sT_\La$ for the $\sigma$-field generated
by the set $\{\om(e):e\in\EE^d\setminus\EE_\La\}$ of states of edges 
having at least one endvertex outside $\Lambda$.

\definition{3.\Thm} Let $0\le p\le 1$ and $q>0$.
A probability measure $\phi$ on $(\Om ,\sF
)$ is called a {\it random-cluster measure\/} with parameters $p$ and
$q$ if
$$
\text{for
all $A\in\sF$ and all finite boxes $\Lambda$,}\q
\phi (A\mid\sT_\Lambda )(\xi)=\phi_{\Lpq}^\xi (A)\q\text{for $\phi$-a.e.\ 
$\xi$.}
$$
The set of such measures is denoted $\sR_{p,q}$.
\label{rcdef}

The condition of this definition amounts to the following.  
Suppose we are given that the configuration off
the finite box $\La$ is that of $\xi$. Then, for almost every $\xi\in\Omega$,
the (conditional) measure
on $\La$ is simply the \rc\ measure with boundary condition $\xi$.
No further generality is gained by replacing the finite box $\Lambda$
by a general finite subset of $\Z^d$.

Some information about the structure of $\sR_{p,q}$, and its relationship
to $\sW_{p,q}$, is provided in [\ci{G93}]. For example, for all $p$, $q$,
$\sR_{p,q}$
is non-empty and convex. We have no proof that $\sW_{p,q}\subseteq\sR_{p,q}$,
but we
state one theorem in this direction.
For $\om\in\Om$, let $I(\om)$ be the number of infinite open
clusters of $\om$. We say that a probability 
measure $\phi$ on $(\Om,\sF)$ has the
$0/1$-{\it infinite-cluster property\/} if $\phi(I\in\{0,1\})=1$.

\proclaim{Theorem 3.\Thm{} [\ci{G93}, \ci{G96b}, \ci{G02}]}  
Let $0\le p \le 1$ and $q>0$. If $\phi\in\ccw$
and $\phi$ has the $0/1$-infinite-cluster property, then $\phi\in\sR_{p,q}$.
\endproclaim
\label{rcexist}

Since, [\ci{BK}], any translation-invariant probability measure 
satisfying a finite-energy property (see 
the discussion after Theorem 3.\refer{tower}) necessarily
has the $0/1$-infinite-cluster property, we have that
all translation-invariant members
of $\ccw$ lie in $\sR_{p,q}$. 
Suppose for the moment that $q\ge 1$. 
By Theorem 3.\refer{thermo}(b), the weak
limits $\phi_{p,q}^b$, $b=0,1$, are translation-invariant, and therefore
they belong to $\sR_{p,q}$. 
It is not difficult to see, by the FKG inequality,  that
$$
\fpq^0 \lest \phi \lest \fpq^1 \qq\text{for all $\phi\in \sR_{p,q}$},
\tag"($\ast$)"
$$
and it follows 
that  $|\sR_{p,q}|=1$ if and only if $\fpq^0=\fpq^1$.
The claim
of ergodicity in Theorem 3.\refer{thermo}(c) is one
consequence of the extremality ($\ast$) of the $\phi_{p,q}^b$ within the
class $\sR_{p,q}$ (see also [\ci{L00}, page 1113]). 

It may be seen by an averaging argument, [\ci{G93}],
that $\ccw$ necessarily contains at least one
translation-invariant measure, for all $p\in[0,1]$ and $q\in(0,\infty)$.
Therefore, $\sR_{p,q}$ is non-empty for all $p$ and $q$.

We note that Theorem 3.\refer{rcexist}, and particularly the
$0/1$-infinite-cluster property,
is linked to the property
of so-called `almost sure quasilocality', a matter
discussed in 
[\ci{PV}].

\subsection{3.6}{The case $q<1$}
The FKG inequality, a keystone of many arguments when $q \ge 1$, is not 
valid when $q<1$. Consequently,
many fundamental questions are unanswered to date, and the theory
of \rc\ models on a finite graph $G=(V,E)$ remains obscure when $q<1$.
The intuition is that certain positive
correlations should be replaced by negative correlations;
however, the theory of negative correlation is more problematic
than that of positive correlation 
(see [\ci{Pem00}]). We return to this point later
in this subsection.

As referred to above, there is an existence proof 
of infinite-volume weak limits and random-cluster measures
for all $q>0$. On the other hand, no constructive proof
is known of the existence of
such measures when $q<1$ \OP. More specifically, the existence of the weak
limits $\lim_{\La\up\Z^d}\phi_{\Lpq}^b$, $b=0,1$, is not known when $q<1$.
The best that can be shown currently is that the two limits exist and 
are equal when $p$ is either sufficiently small or sufficiently large,
[\ci{G02}].
This may be achieved by comparison with percolation
models having different values of $p$, very much as in [\ci{G93}]
(the claim for small $p$ may also be shown by the arguments
of [\ci{FFG1}, \ci{FFG2}]).

The theory of percolation gives a clue to a possible way forward.
When $q=1$, the FKG inequality is complemented by the so-called
`disjoint-occurrence' (or `BK') inequality. This latter inequality 
is said to be valid for a measure $\mu$ if
$\mu(A\circ B) \le \mu(A)\mu(B)$ for all increasing events $A$, $B$,
where $A\circ B$ is the event that $A$ and $B$ occur disjointly
(see [\ci{vdBK}, \ci{G99}] for a discussion of this and the more general
`Reimer inequality' [\ci{Rei}]). 
The disjoint-occurrence inequality has been established
for classes of measures which are
only slightly more general than product measures, and it is 
an interesting open question whether it is valid for a wider
class of measures of importance \OP. It has been asked
whether the disjoint-occurrence inequality could be valid for
\rc\ measures with $q<1$ \OP. A positive answer would 
aid progress substantially towards an understanding of limit \rc\ measures.

We illustrate this discussion about disjoint-occurrence
with the following test question \OP: is
it generally the case that  the \rc\ measure $\fpq$ on $G$ satisfies
$$
\fpq(\text{edges $e$ and $f$ are open})
\le\fpq(\text{$e$ is open})\fpq(\text{$f$ is open})
\tag$\ast$
$$
for $e\ne f$ and $q < 1$? (See [\ci{Pem00}].)
This equation would be a very special instance of the disjoint-occurrence
inequality. A further restriction arises if we take the limit
as $q\downarrow 0$; recall the discussion of Section 2.4.  
This leads to certain open questions
of a purely graph-theoretic type, which combinatorial
theorists might elevate to the status of conjectures. The
first such question is the following.
Let $K(e_1,e_2,\dots)$ be the number of subsets $F$ of the edge set $E$, containing
$e_1,e_2,\dots$, such that the graph $(V,F)$ is connected.
Is it the case that \OP
$$
K(e,f)K(\es) \le K(e)K(f)\quad\text{if } e\ne f?
\tag$\ast\ast$
$$
(See [\ci{Kahn}].)
In the second such question, we ask if the same inequality is valid with
$K(e_1,e_2,\dots)$ redefined as the number of
subsets $F$ containing $e_1,e_2,\dots$
such that $(V,F)$ is a forest \OP. These two questions are
dual to one another in the sense that the first holds for a planar
graph $G$ if and only if the second holds for its planar dual.
Explicit computations have confirmed the forest conjecture for all graphs
$G$ having nine or fewer vertices, [\ci{GW}].

In the `intermediate regime', with $K(e_1,e_2,\dots)$ redefined
as the number of spanning trees (that is, connected forests) of $G$
containing $e_1,e_2,\dots$,
the corresponding inequality is indeed valid. An extra ingredient
in this case is the link to electrical networks, and 
particularly the variational
principle known as the Thomson or Dirichlet principle (see [\ci{DoyS}]).
Further results and references are provided in [\ci{BLPS01}]. 
Substantially more is known
for spanning trees, namely a general result concerning the `negative
association' of the uniform measure on the set of spanning trees
of $G$, [\ci{FM92}].

We note a more general version of conjecture ($\ast\ast$), namely
$$
K_\a(e,f)K_\a(\es) \le K_\a(e)K_\a(f)\quad\text{for } e\ne f,\ 0<\a<\infty,
$$
where
$$
K_\a(e_1,e_2,\dots)=\sum\Sb F\subseteq E\\ F\supseteq\{e_1,e_2,\dots\}\\
   (V,F) \text{ connected}\endSb \a^{|F|}.
$$
This is equivalent to ($\ast$) in the limit as $q\downarrow 0$,
where $\a=p/(1-p)$.

By other means one may establish
a certain non-trivial monotonicity when $q<1$, but by a
more complicated reasoning than before involving a
property of convexity of the logarithm of the partition function. Namely,
the mean number of open edges is non-decreasing
in $p$, for $0<q<\infty$,
[\ci{G93}].

\endsection

\section{4. Phase transition, the big picture}
Phase transition in a Potts model corresponds to the creation
of an infinite open cluster in the corresponding \rc\ model.
There are rich predictions concerning the nature of
such a phase transition, but these have been  
proved only in part.
This section is a summary of the expected properties of the phase diagram for
different dimensions $d$ and cluster-weighting factors $q$.
The corresponding rigorous theory is described in Sections 5 and 6.

\subsection{4.1}{Infinite open clusters}
{\it We assume henceforth that\/}
$q\geq 1$, and we concentrate here on the extremal random-cluster measures
$\phi_{p,q}^0$ and $\phi_{p,q}^1$. 
The phase transition of a random-cluster measure is
marked by the onset of an infinite
open cluster.
We write $\{0\lra\infty\}$ for
the event that the origin is the endvertex of some infinite open path,
and we define the $\phi_{p,q}^b$ {\it percolation probability\/} by
$$
\t^b (p,q)=\phi_{p,q}^b (0\lra\infty ),\qq b=0,1.
$$
It is almost immediate by a stochastic-ordering argument
that $\t^b (p ,q)$ is non-decreasing in $p$, and therefore 
$$
\t^b(p,q) \cases =0 &\text{if } p < \pc^b(q),\\
>0 &\text{if } p > \pc^b(q),
\endcases \qq b=0,1,
$$
for critical points $\pc^b(q)$  given by
$$
\pc^b(q)=\sup\{ p:\t^b(p,q)=0\},\qq b=0,1.
$$
It is an easy exercise to show that the number $I$ of infinite
open clusters satisfies:
$$
\fpq^b(I\ge 1)= \cases 0 &\text{if } \t^b(p,q)=0,\\
1 &\text{if }\t^b(p,q)>0.
\endcases
$$
We shall see in Section 5.2 that any infinite 
open cluster is $\fpq^b$-a.s.\ unique
whenever it exists.

We have by Theorem 3.\refer{countable} that $\phi_{p,q}^0=\phi_{p,q}^1$ for
almost every $p$, whence $\t^0(p,q)=\t^1(p,q)$ for almost every $p$,
and therefore $\pc^0(q)=\pc^1(q)$. Henceforth we use the abbreviated
notation
$\pc(q)=\pc^0(q)=\pc^1(q)$,
and we refer to $\pc(q)$ as the {\it critical point\/} of the
corresponding random-cluster measures. The non-triviality of
$\pc(q)$ may be proved by comparisons of \rc\ measures with product measures
via Theorem 3.\refer{dom}. Recall the fact, [\ci{G99}, Chapter 1],
that $0<\pc(1)<1$ if $d \ge 2$.

\proclaim{Theorem 4.\Thm{} [\ci{ACCN}]} 
We have for $q\ge 1$ that
$$
\pc(1) \le \pc(q) \le \frac{q\pc(1)}{1+(q-1)\pc(1)}.
$$
\endproclaim
\label{nontriv}

When $q$ is an integer satisfying $q\ge 2$, the phase transition of 
the random-cluster model corresponds in a special way to that
of the Potts model with the same value of $q$.
An indicator of phase transition in the Potts model is the
`magnetisation', defined as follows. Consider a Potts measure
$\pi_\La^1$ on $\La$ having parameters $\b$, $J$, $q$,
and with `1' boundary conditions, which is
to say that all vertices on the boundary $\partial\Lambda$ are
constrained to have spin value 1. Let $\tau_\Lambda =\pi_\Lambda^1 (\s_0
=1)-q^{-1}$, a quantity which represents the net effect 
of this boundary condition on the spin
at the origin. The corresponding random-cluster measure
$\phi_\La^1$ has parameters $p=1-e^{-\b J}$ and $q$, and 
has wired boundary condition.
We apply Theorem 2.\refer{corrconn} to the graph obtained from
$\Lambda$ by identifying all vertices in $\partial\Lambda$,
and we find that
$$
\tau_\Lambda =(1-q^{-1})\phi_\Lambda^1(0\lra\partial\Lambda ).
$$
The limit function
$\tau\,
=\lim_{\Lambda\up\Z^d} \tau_\Lambda$ is called the magnetisation,
it is a non-decreasing function of $\b J$ and satisfies
$$
\tau\  \cases =0 &\text{if  $\b J$ is small},\\
>0 &\text{otherwise}.
\endcases
$$
It
is not hard to show, [\ci{ACCN}], that $\phi_\Lambda^1(0\lra\partial\Lambda )
\to \phi^1(0\lra\infty)$ as $\La\up\Z^d$, whence $\tau=(1-q^{-1})\t^1(p,q)$
where $p=1-e^{-\b J}$. 
Therefore there
is long-range order in the Potts model (that is, $\tau>0$) if and only if
the origin lies in an infinite open cluster with strictly
positive $\phi_{p,q}^1$-probability. In particular, 
$\pc(q)=1-e^{-\beta_{\text{c}}J}$ where $\beta_{\text{c}}$ is the critical
value of $\beta$ for the Potts model in question.

\subsection{4.2}{First- and second-order phase transition}
There is a rich physical theory of phase transitions in percolation,
Ising, and Potts models, some of which has been made rigorous in the 
context of the random-cluster model. There follows a broad
sketch of the big picture, a full rigorous verification of
which is far from complete. Rigorous 
mathematical progress is described in
Section 5.

\ssubsection{I}{The subcritical phase, $p<\pc(q)$.}
It is standard,  [\ci{ACCN}], that
$$
\fpq^0 =\fpq^1\qq\hbox{if $\t^1 (p,q)=0$},
$$
implying that there exists a unique random-cluster measure whenever $\t^1
(p,q)=0$. In particular,
$|\sW_{p,q}|=|\sR_{p,q}|=1$ when $0\le p < \pc(q)$.
Assume for the moment that $p<\pc(q)$, and denote 
the unique \rc\ measure by $\fpq$.
By the definition of the critical point, all open clusters are $\fpq$-a.s.\
finite. It is believed that they have a tail which decays exponentially,
in that there  exist $\g=\g(p,q) > 0$ and $\eta=\eta(p,q)>0$  such that
$$
\fpq(|C|=n) = e^{-\g n(1+\o(1))},\q \fpq(\rad(C)=n) =  e^{-\eta n(1+\o(1))},
\qq\text{as } n\to\infty,
$$
where $C$ denotes the open cluster containing the origin,
and its {\it radius\/} 
$\rad(C)$ is defined as $\sup\{\|x\|:x\in C\}$.
Such exponential decay would be the starting point for a complete exploration
of the subcritical phase. More detailed asymptotics should then emerge,
including the Ornstein--Zernike decay of the connectivity functions:
$$
\fpq(0\lra e_n)\sim \frac{c(p,q)}{n^{(d-1)/2}}e^{-n/\xi(p,q)}
  \qq\text{as } n\to\infty,
$$
where $e_n=(n,0,0,\dots,0)$ and $\xi(p,q)$ is termed the `correlation length'.

\ssubsection{II}{The supercritical phase, $p>\pc(q)$.}
This phase is characterised by the existence of one or more infinite
open clusters (exactly one, in fact, for translation-invariant
measures at least, see Section 5.2).
It is believed that, as in the subcritical phase,
we have that $\fpq^0=\fpq^1$ when
$p>\pc(q)$; this remains unproven in general. Thus the first
main problem is to prove that there is a unique random-cluster measure
when $p>\pc(q)$. 

The theory of percolation, [\ci{G99}],
suggests a route towards understanding the geometry of the
supercritical phase, namely by developing a rigorous block renormalisation
argument.  This should permit the use of
theory developed when $p$ is close to 1 in order to
understand the model when  $p$ is close to $\pc(q)$.
In particular, one expects an exponential estimate for the
decay of the probabilities of long-range connections within
finite open
clusters, and a Wulff construction for the shape of  
such clusters.

\ssubsection{III}{Near the critical point, $p\simeq\pc(q)$.}
The main open problem is to understand the way in which the
nature of the phase transition depends on the value of $q$.
It is believed that the transition is continuous and governed
by critical exponents and scaling theory when $q$ is small,
and is discontinuous when $q$ is large. Presumably there exists   
a threshold for $q$ which separates 
the so-called `second-order' (or continuous)
transition from the so-called `first-order' (or discontinuous) transition.
More specifically, it is believed that there
exists $Q=Q(d)$ satisfying
$$
Q(d)=\cases 4 &\text{if } d = 2,\\
2 &\text{if } d \ge 6,
\endcases
$$
such that the following hold.
\ii{(i)} {\it Assume that  $q<Q$.} 
\iii{--}For any $p$, there exists a 
unique random-cluster measure, denoted $\fpq$.
In particular $\phi_{\pc(q),q}^0=\phi_{\pc(q),q}^1$. 
\iii{--}
$\t(p,q)=\fpq(0\lra\infty)$ is a continuous function of $p$. 
There is no percolation at the critical point, in the sense that
$\t(\pc(q),q)=0$. 
\iii{--}
The edge-density $h(p,q)=\fpq(e \text{ is open})$, viewed as a function
of $p$, 
is continuous at the critical point $p=\pc(q)$. [The letter $e$ denotes
a typical edge of the lattice.] 
\iii{--}
These functions and others have power-law
singularities at $\pc(q)$, and the associated critical exponents 
satisfy the scaling relations (see [\ci{G99}, Chapter 9]). 
\iii{--} When $d$ is large (how large depends on the value of $q$), 
these critical
exponents take on their `mean-field' values, and
depend no further on the value of $d$.
\iii{--} There is no `mass gap', in the sense that the correlation length
$\xi(p,q)$ satisfies $\lim_{p\uparrow\pc(q)}\xi(p,q)=\infty$.
\iii{--} Universality reigns, in that the critical exponents depend
on the number $d$ of dimensions but not on the choice of lattice.
For example, the exponents associated with the square lattice are 
expected to be the same as
those for the triangular lattice.
\iii{--} Assume $d=2$ and $1\le q<4$.
The process with $p=\pc(q)$ converges as the lattice spacing
shrinks to zero, the limit process when suitably defined being
a stochastic L\"owner evolution \SLE$_\kappa$ having parameter $\kappa$
satisfying $\cos(4\pi/\kappa)=-\frac12\sqrt q$, $\kappa\in(4,8)$
(see Section 6.4 and [\ci{RS01}]).

\ii{(ii)} {\it Assume that $q>Q$.}
\iii{--} There exists a 
unique random-cluster measure if and only if $p\ne \pc(q)$.
When $d=2$ and $p=\pc(q)$, there are exactly two extremal members 
of $\sR_{p,q}$, namely the free and the wired measures
$\fpq^b$, $b=0,1$. When $d \ge 3$
and $p=\pc(q)$ there exist other extremal members of $\sR_{p,q}$ including
a variety of non-translation-invariant measures.
\iii{--} We have that $\t^0(\pc(q),q)=0$ but $\t^1(\pc(q),q)>0$.
\iii{--} The edge-density $h(p,q)$ is a discontinuous 
function of $p$ at the critical point
$\pc(q)$.
\iii{--} There is a `mass gap' in the sense that the correlation length
$\xi(p,q)$ satisfies $\lim_{p\uparrow\pc(q)}\xi(p,q)<\infty$.

\endsection

\section{5. General results in $d$ ($\ge 2$) dimensions}
The properties of the \rc\ model depend pivotally on whether the process
is subcritical ($p<\pc(q)$), supercritical ($p>\pc(q)$),
or critical ($p\simeq\pc(q)$). We consider
these situations in turn, in each case  
identifying major results and open problems. 
There is a bulk of information
available for certain values of $q$, 
namely when $q=1,2$ and $q$ is sufficiently large. In addition,
the case $d=2$ is special, and we shall return to this in Section 6.
We assume throughout this section that $q\ge 1$.

Little is known in general
about the numerical values of $\pc(q)$. For example, it is known that $\pc(q)$
is Lipschitz-continuous and strictly increasing when $d \ge 2$, [\ci{G95a}],
and there is a striking conjecture \OP\ that $\pc(q)=\sqrt q/(1+\sqrt q)$
when $d=2$ (see Section 6.2). Some concrete inequalities involving
the $\pc(q)$ are implied by the comparison inequalities of
Theorem 3.\refer{dom}.
 
\subsection{5.1}{The subcritical phase, $p<\pc(q)$}
As remarked in Section 4.2, there is a unique \rc\ measure when $p<\pc(q)$,
and we shall denote this by $\fpq$.

The key theorem for understanding the subcritical phase of percolation
states that long-range connections have exponentially decaying
probabilities. Such a result is believed to hold for all
random-cluster models with $q\geq 1$, but no proof has been
found \OP\ which is valid for all $q\ge 1$ and all $p<\pc(q)$. 
The full result is known only when $q=1$, $q=2$, or $q$ is
sufficiently large, and the three sets of arguments for these cases
are somewhat different from one another. As for results valid for all
$q$ ($\geq 1$), the best that is currently known
is that the connectivity function decays  exponentially for
sufficiently small $p$ (this follows  by
Theorem 3.\refer{dom} and the corresponding
$q=1$ result), and that it decays {\it exponentially\/}
whenever
it decays at a sufficient {\it polynomial\/} rate.
We describe the last result next.

As a preliminary we introduce another definition of a critical point.
Let $B(n)$ be the cube $[-n,n]^d$. We write
$$
Y(p,q)=\limsup_{n\to\infty}
\Bigl\{ n^{d-1}\fpq\bigl( 0\lra\pd B(n)\bigr)\Bigr\}
$$
and
$\pcg(q)=\sup\bigl\{ p:Y(p,q)<\infty\bigr\}$.
Evidently $0<\pcg(q)\leq\pc (q)$, and it is believed that 
$\pcg(q)=\pc(q)$ for all $q\ge 1$ \OP.

\proclaim{Theorem 5.\Thm{} [\ci{GP96}]} Let $q\ge 1$, $d\ge 2$,
and $0\le p<\pcg(q)$.
There exists $\gamma=\gamma(p,q)$
satisfying $\gamma > 0$ such that
$$
\fpq( 0\lra\pd B(n))\leq e^{ -\gamma n}
\qq\text{for all large } n.
$$
\endproclaim
\label{expdecay}

The spirit of the theorem is close to that of Hammersley [\ci{H57a}]
and Simon--Lieb [\ci{Li80}, \ci{Si80}] who proved exponential estimates
when $q = 1,2$ subject to a hypothesis of finite susceptibility (that is,
under the hypothesis that
$\sum_x\fpq (0\lra x)<\infty$). The latter assumption is slightly
stronger than the assumption of the above theorem when $d=2$.

Connectivity functions are expected to decay exponentially
with a correction term of power order. More specifically,
it is expected as reported in Section 4.2 that 
$$
\fpq(0\lra x) \sim 
\frac c{|x|^{(d-1)/2}}\exp(-|x|/\xi)\qquad\text{as } |x|\to\oo,
$$
for constants $c(p,q)$ and $\xi(p,q)$, and for some
suitable norm $|\cdot|$ on $\Z^d$. Such `Ornstein--Zernike' decay is
a characteristic of many systems in their disordered
phases. No proof is known \OP, except in the special
cases when $q=1$ and $q=2$, [\ci{CCC}, \ci{CIV}]. In [\ci{Al01b}]
may be found a weaker result which bounds the fluctuations
by a power-law when $d=2$, under the assumption that
the function does indeed decay exponentially.

\subsection{5.2}{The supercritical phase, $p>\pc(q)$}
We assume as usual that $q\ge 1$, and we begin with a discussion
of the number of infinite clusters. For $\om\in\Om$, let
$I(\om)$ be the number of infinite open clusters.
\comment
It is an important open problem to prove that $\fpq^0=\fpq^1$ when $p>\pc(q)$,
or equivalently that there exists a unique random-cluster measure
throughout the phase \OP. The best that has been achieved to date is
uniqueness when $p$ is sufficiently close to 1, and this was done, [\ci{G93}],
by making a comparison with percolation at high density. 
\endcomment
Suppose that $\phi_{p,q}$ is a translation-invariant member of $\sR_{p,q}$.
If in addition $\phi_{p,q}$ is ergodic,
then, by a well known theorem of
Burton and Keane [\ci{BK}],
$$
\text{either}\q \phi_{p,q}(I=0)=1 \q\text{or}\q \phi_{p,q}(I=1)=1;
$$
that is to say, the infinite open
cluster is almost surely unique whenever it exists.
It is noted in [\ci{BK}] that methods of ergodic decomposition
enable the extension of such results to translation-invariant
measures which are not necessarily ergodic. That is, under the assumption
of translation-invariance alone,
$$
\phi_{p,q}(I\in\{0,1\})=1,
$$
which is to say that translation-invariant \rc\ measures have the
$0/1$-infinite-cluster property.
A further comment on the use of ergodic decomposition in this context is
to be found in [\ci{BK2}].

\comment
It is rather harder to prove so-called
`simultaneous uniqueness'. One needs first a suitable coupling of \rc\ measures.
It was shown in [\ci{G93}]  that there exists a Markov process 
$\zeta=(\zeta_t:t\ge 0)$ on
the state space $[0,1]^{\EE^d}$ with the following properties.
For $0\le p\le 1$ and $\sigma\in[0,1]^{\EE^d}$, we define
$$
\sigma^{p,-}(e) = 1_{\{\zeta_t(e)>1-p\}},\q
\sigma^{p,+}(e) = 1_{\{\zeta_t(e)\ge 1-p\}},
\qq e\in \EE^d,
$$
where $1_A$ denotes the indicator function of
an event $A$. Then
$\zeta_t^{p,-}$ and $\zeta_t^{p,+}$
are Markov processes
which are reversible with respect to the
infinite-volume free and
wired random-cluster measures, $\fpq^0$ and $\fpq^1$, respectively. 
The process $\zeta_t$ has a unique invariant measure $\mu$, and 
the question of simultaneous uniqueness amounts 
in this context to asking 
whether or not
$$
\mu\Bigl(\text{for all $p>\pc(q)$, } \sigma^{p,-} \text{  and } \sigma^{p,+}
\text{ have unique infinite open clusters}\Bigr)=1,
$$
where $\sigma$ is to be regarded as an element chosen according to $\mu$.
The main theorem of [\ci{Al95}] may be applied
in this context to deduce that this is indeed true. 
A little extra work is necessary to verify the finite-energy
condition of [\ci{Al95}], and this has been done by [\ci{Al02b}, \ci{G02}].
\endcomment

In two dimensions, the supercritical process is best studied via
the subcritical process which arises as its graphical dual 
(see Section 6).
There are two general approaches to the
supercritical phase in a general number $d$ ($\ge 3$) of dimensions.
The less powerful is to derive results for large $p$ by comparison
with percolation, the theory of which is
relatively complete. Without an extra ingredient, such an approach
will not reveal the structure of the supercritical phase all the way down to
the critical value $\pc(q)$.  As an example, we present one theorem
concerning the uniqueness of \rc\ measures.

\proclaim{Theorem 5.\Thm{} [\ci{G93}]} If $d\ge 2$ and $q\ge 1$, there exists
$p'=p'(d,q)<1$ such that $\phi_{p,q}^0=\phi_{p,q}^1$
whenever $p>p'$.
\endproclaim
\label{unique}

It is an important open problem to prove that $\fpq^0=\fpq^1$ 
for all $p>\pc(q)$,
or equivalently that there exists a unique random-cluster measure
throughout the phase \OP. 

A more powerful approach, sometimes used in conjunction with
the comparison argument summarised above, is the `block argument'
laid out in [\ci{CP}, \ci{Pisz96}].
One may think of block arguments as a form
of rigorous renormalisation.
One divides
space into blocks, constructs events of an appropriate 
nature on such blocks, having large probabilities,
and then allows these events to combine across space. 
There have 
been substantial successes using this technique, 
of which the most striking is the resolution, 
subject to certain side conditions,
of the so-called Wulff construction for the asymptotic shape of large
Ising droplets.

Rather than discussing the physical background of the Wulff
construction, we mention instead its impact on \rc\ models.
Let $B(n)=[-n,n]^d$, and consider the wired \rc\ measure
$\phi_{B(n),p,q}^1$  with $p>\pc(q)$.
The larger is an open cluster, the more likely it is
to be joined to the boundary $\partial B(n)$.
Suppose that we condition on the event that
there exists in $B(n)$ an open cluster $C$ which does
not touch $\partial B(n)$ and which has volume of the order
of the volume $n^d$ of the box. What can be said about the
shape of $C$? Since $p>\pc(q)$, there is little cost in having 
large {\it volume\/},
and the price is spent around its {\it boundary\/}.
Indeed, the price may be expressed as a surface integral
of an appropriate function termed `surface tension'.
This `surface tension' may be specified as the
exponential rate of decay of a certain probability.
The Wulff prediction for the shape of $C$ is that, when re-scaled
in the limit of large $n$, its shape converges to the solution
of a certain variational problem, that is, the limit shape
is obtained by minimising a certain surface integral subject
to a condition on its volume.

No proof of this general picture for \rc\ models has appeared in
the literature, although it is believed that the methods
of [\ci{CP}, \ci{CP01}, \ci{Pisz96}] enable such a proof. 
The authors of [\ci{CP}]
have instead concentrated on using \rc\ technology to
solve the corresponding question for the asymptotic shape of
large droplets in the Ising model.
The outcome is an important `large deviation' theorem
which utilises block arguments and yields a full
solution to the Ising problem whenever the corresponding random-cluster model 
(which has $q=2$) has parameter $p$ satisfying
$p>\wpc(2)$ and $\phi^0_{p,2} = \phi^1_{p,2}$.
Here, $\wpc(2)$ is the limit of a certain decreasing
sequence of critical points defined on slabs in $\Z^d$,
and is conjectured \OP\ to be equal to the critical point $\pc(2)$.
[Closely related results have been obtained in [\ci{Bod99}].
Fluctuations in droplet shape for {\it two-dimensional\/} \rc\ models
have been studied in [\ci{Al02}, \ci{Al02a}].]

The `slab critical point' $\wpc(q)$ may be defined for any \rc\ model
as follows.
Fix $q\ge 1$, and let $d\ge 3$. 
Let $S(n,L)=[-n,n]^{d-1}\times[-L,L]$.
Let $\psi_{p,q}^{n,L}$ be the \rc\ measure on $S(n,L)$ with parameters
$p$, $q$ (and with free boundary conditions).
We denote by $\Pi(p,L)$ the property that:
$$
\text{there exists $\a>0$
such that, for all } x\in S(n,L) \text{ and all } n,
\psi_{p,q}^{n,L}(0\lra x) > \a.
$$
It is not hard to see that $\Pi(p,L)\Rightarrow \Pi(p',L')$
if $p\le p'$ and $L\le L'$. It is thus natural to define the quantities
$$
\wpc(q,L)=\inf\{p: \Pi(p,L)\text{ occurs}\},\qq
\wpc(q)=\lim_{L\to\infty} \wpc(q,L),
$$
and it is clear that
$\wpc(q)\ge \pc(q)$.

\proclaim{Conjecture 5.\Thm{} [\ci{Pisz96}]} Let $q\ge 1$ and $d \ge 3$.
We have that $\wpc(q)=\pc(q)$.
\endproclaim
\label{pchat}

Subject to a verification of this conjecture, and of a positive
answer
to the question of the uniqueness of \rc\ measures when $p>\pc(q)$,
the block arguments of [\ci{CP}, \ci{Pisz96}] may be expected to result in
a fairly complete picture of the supercritical phase of \rc\ models with
$q\ge 1$; see [\ci{CP01}] also.

The case $q=1$ is special, percolation enjoys a spatial independence
not shared with general \rc\ models. This additional property has
been used  in the  formulation of a type of `dynamic renormalisation',
which has in turn yielded a proof that
$\wpc(1)=\pc(1)$ for percolation in three
or more dimensions,
[\ci{G99}, Chapter 7, \ci{GrM}].
Such arguments do not to date have a \rc\ counterpart.

As a further application of a block argument
we note the following bound, [\ci{Pisz96}], 
for the tail of the size of
the open cluster $C$ at the origin,
$$
\fpq^b(|C|=n) \le \exp\bigl(-\a n^{(d-1)/d}\bigr)\qq\text{for all } n,
$$
for some $\a=\a(p,q)>0$, and
valid for $d\ge 3$, $b=0,1$, and $p$ sufficiently close to 1.
The complementary inequality
$$
\fpq^b(|C|=n) \ge \exp\bigl(-\a' n^{(d-1)/d}\bigr)\qq\text{for all } n,
$$
may be obtained for large $p$ as done in the
case of percolation, [\ci{G99}, Section 8.6].

\subsection{5.3}{Near the critical point, $p\simeq\pc(q)$}
Surprisingly little is known about \rc\ measures near the critical
point, except in the cases $q=1,2$ and $q$ large. In each such case,
there are special arguments which are apparently not suitable
for generalisation.
We summarise such results as follows.

\ssubsection{I}{Percolation, $q=1$.}
There is a full theory of the subcritical and supercritical phases
of percolation, [\ci{G99}]. The behaviour when $p\simeq \pc(1)$ has been the
subject of deep study, and many beautiful results are known. 
Nevertheless, the picture
is incomplete. For example, it is believed but not proved
that $\t(\pc(1),1)=0$ for all $d\ge 2$, but this is known only when $d=2$ (because of
special properties of two dimensions explored for $\LL^2$ in Section 6)
and when $d$ is large ($d\ge 19$ suffices) using a method termed the
`lace expansion'. The lace expansion explains also the values of 
some critical exponents
when $d$ is large; see, for example, [\ci{HS}, \ci{HS00}].

Great progress has been made in recent years towards understanding
the phase transition when $d=2$.  The idea is to work at
the critical point $p=\pc(1)$,
and to observe the process over an increasing sequence of regions
of $\Z^2$. It is believed that the process, re-scaled as
the regions become larger, converges in a certain manner to
a stochastic process generated in a prescribed way 
by a differential
equation, known as a L\"owner equation, which is driven 
in a certain way by a Brownian motion.
Stochastic processes which arise in this way
have been termed {\it stochastic L\"owner evolutions\/} by
Schramm, [\ci{Sch00}], and denoted \SLE$_\kappa$,
where $\kappa$ is the variance parameter of the Brownian motion.
It is believed that the space of stochastic L\"owner evolutions is a 
canonical family of processes which arise as scaling limits
of discrete processes such as critical percolation,
critical \rc\ models with $q\le 4$, self-avoiding walks,  
loop-erased random walk, and
uniform spanning trees. Full proofs are not yet known \OP.
We expand on this very important development in Section 6.4

\ssubsection{II}{Ising model, $q=2$.}
Integer values of $q$ are special, and the value $q=2$ 
particularly so because of certain transformations
which permit the passage to a model which might be termed a `Poisson graph'.
Let $G=(V,E)$ be a finite graph and let $0<\l<\oo$. Suppose that $\pi=
\{\pi(e):e\in E\}$ is a family of independent random variables
each having the Poisson distribution with parameter $\l$.
We now construct a random graph $G_\pi=(V,E_\pi)$ having vertex
set $V$ and, for each $e\in E$, having exactly $\pi(e)$ edges
in parallel joining the endvertices of the edge $e$
[the original edge $e$ is itself removed].
We call $G_\pi$ a {\it Poisson graph with intensity $\l$}, and write
$\PP_\l$ and $\EE_\l$ for the appropriate probability measure
and expectation operator.

We introduce next the concept of a flow on an oriented graph.
Let $q\in\{2,3,\dots\}$ and let $G'=(V',E')$ be a finite oriented graph. 
Let $f:E'\to
\{0,1,2,\dots,q-1\}$. For $x\in V'$, the {\it total flow into $x$\/}
is the sum of $\pm f(e')$ over all edges $e'$ incident to $x$, with
$+1$ when $e'$ is oriented towards $x$ and $-1$ otherwise.
The function $f$ is called a {\it mod-$q$ flow\/} if the total flow
into $x$ is zero (modulo $q$) for all $x\in V'$.
The mod-$q$ flow $f$ is called {\it non-zero\/} if $f(e')\ne 0$
for every $e'\in E'$. We write $F_q(G')$ for the number
of non-zero mod-$q$ flows on $G'$. It is a remarkable fact, [\ci{Tut}], that
$F_q(G')$ does not depend on the orientations of 
edges in $E'$, and thus
one may define $F_q(G')$ unambiguously for any {\it unoriented\/} graph $G'$.

We return now to the Poisson graph $G_\pi$. For $x,y\in V$, $x\ne y$,
we denote by $G_\pi^{x,y}$ the graph obtained from $G_\pi$ by adding an edge
with endvertices $x$, $y$. [If $x$ and $y$ are already 
adjacent in $G_\pi$, we add exactly
one further edge between them.] 
Connection probabilities and flows are related by the following theorem,
which may be proved using properties of Tutte polynomials 
(see [\ci{Tut}] and Section 2.5).

Let $G=(V,E)$
be a finite graph, and write $\phi_{G,p,q}$ for the \rc\ measure on $G$ with
parameters $p$, $q$.

\proclaim{Theorem 5.\Thm{} [\ci{G91}, \ci{G02}]} 
Let $q \in\{2,3,\dots\}$ and $0\le p=1-e^{-\lambda q}< 1$. 
We have that
$$
(q-1)\phi_{G,p,q}(x\lra y) = \frac{\EE_\l(F_q(G_\pi^{x,y}))}{\EE_\l(F_q(G_\pi))}
\qq\text{for all } x,y\in V,\ x\ne y.
$$
\endproclaim
\label{flows}

This formula takes on an especially simple form when $q=2$, since
non-zero mod-2 flows necessarily take only the value 1. It follows that,
for any graph $G'$, $\EE_\l(F_2(G'))$ equals the $\PP_\l$-probability
that the degree of every vertex of $G'$ is even, [\ci{A82}].
Observations of this sort have led when $q=2$ to the so-called 
`random-current' expansion for Ising models, thereby
after some work [\ci{A82}, \ci{ABF}, \ci{AF}] 
leading to proofs amongst other things
of the following, expressed here in the language of \rc\ measures.

\ii{(i)} When $q=2$ and $p<\pc(q)$, we have 
exponential decay of the radius distribution,
$$
\phi_{p,2}(\rad(C)=n)\le e^{-\eta n}\qq\text{for all } n,
$$
where $\eta=\eta(p)>0$; exponential decay of the
two-point connectivity function follows.

\ii{(ii)} When $q=2$ and $d\ne 3$, there is a unique
\rc\ measure $\phi_{p,2}$ for all $p$,
in that $|\sR_{p,q}|=1$.

\ii{(iii)} The phase transition is continuous 
when $q=2$ and $d\ne 3$. In particular,
$\t^0(\pc(2),2)=\t^1(\pc(2),2)=0$, and the edge-density 
$h(p,2)=\phi_{p,2}(e\text{ is open})$
is a continuous function of $p$ at the critical point $\pc(2)$.

\ii{(iv)} When $d\ge 4$, some (at least) critical exponents take their
mean-field values, and depend no further on the value of $d$.

Note that the nature of the phase transition in three dimensions 
remains curiously undecided \OP.

\ssubsection{III}{The case of large $q$.}
It is not known whether the phase transition is continuous for all small $q$
\OP.
The situation for large $q$ is much better understood owing to a method
known as Pirogov--Sinai theory [\ci{PS1}, \ci{PS2}] 
which may be adapted in a convenient
manner to \rc\ measures. The required computation, which may be found in
[\ci{LMMRS}], has its roots in an earlier paper [\ci{Kot-S}]
dealing with Potts models. A feature of such arguments is 
that they are valid `all the way to the critical point'
(rather than for `small $p$' or `large $p$' only), 
so long as $q$ is sufficiently large. 
One obtains thereby a variety of conclusions including
the following.

\ii{(i)} The edge-densities $h^b(p,q)=\fpq^b(e\text{ is open})$, $b=0,1$, 
are discontinuous functions
of $p$ at the critical point.
\ii{(ii)} The percolation probabilities satisfy
$\t^0(\pc(q),q)=0$, $\t^1(\pc(q),q)>0$.
\ii{(iii)} There is a multiplicity of \rc\ measures
when $p=\pc(q)$, in that $\phi_{\pc(q),q}^0 \ne \phi_{\pc(q),q}^1$.
\ii{(iv)} If $p<\pc(q)$, there is exponential decay and a mass gap, in that
the unique \rc\ measure satisfies
$$
\fpq(0\lra e_n) = e^{-(1+\o(1))n/\xi}\qq\text{as } n\to\oo,
$$
where $e_n=(n,0,0,\dots,0)$ and
the correlation length $\xi=\xi(p,q)$ is such that
$\lim_{p\uparrow\pc(q)}\xi(p,q)=\psi(q)<\oo$.
\ii{(v)} If $d=3$ and $p=\pc(q)$, there exists a non-translation-invariant
random-cluster measure, [\ci{CK}, \ci{MMRS}].

It is not especially fruitful to seek numerical 
estimates on the required size $Q(d)$ of $q$ for
the above conclusions to be valid. Such estimates may be computed,
but turn out to be fairly distant from those anticipated,
namely $Q(2)=4$, $Q(d)=2$ for $d \ge 6$.

The proofs of the above facts are rather complicated and will not
be explained here. Proofs are much easier
and not entirely dissimilar when $d=2$, and a very 
short sketch of such a proof
is provided in Section 6.3.
\endsection

\section{6. In two dimensions}
The duality theory of planar graphs provides a 
technique for studying \rc\ models in two dimensions. 
We shall see in Section 6.1
that, for a dual pair $(G,G^\tD)$ of planar graphs,
the measures $\phi_{G,p,q}$ and $\phi_{G^\tD,p^\tD,q}$
are dual measures in a certain geometrical sense, 
where $p$, $p^\tD$ are related by
${p^\tD}/(1-p^\tD) = {q(1-p)}/{p}$.
Such a duality permits an analysis by which many results
for $\LL^2$ may be derived.  Of particular interest is the value
of $p$ for which $p=p^\tD$.  This `self-dual point'
is easily found to be $p=\psd q$ where
$$
\psd q=\frac{\sqrt q}{1+\sqrt q},
$$
and it is conjectured that $\pc(q)=\psd q$ for $q\ge 1$.

\subsection{6.1}{Graphical duality}
Let $G=(V,E)$ be a simple planar graph
imbedded in $\R^2$.  We obtain its dual graph $G^\tD
=(V^\tD ,E^\tD )$ as follows 
(the roman letter `d' denotes `dual' rather than number
of dimensions). 
We place a dual vertex within
each face of $G$, including the infinite face of $G$ if $G$ is finite.
For each $e\in E$ we place a dual edge $e^\tD=\langle x^\tD,y^\tD\rangle$
joining the two dual vertices lying in the two faces of $G$ abutting $e$;
if these two faces are the same, then $x^\tD=y^\tD$ and $e^\tD$ is a loop.
Thus $E^\tD$ is in one--one correspondence to $E$.
It is easy to see that the dual of $\LL^2$ is isomorphic to $\LL^2$. What
is the relevance of graphical duality to random-cluster measures on $G$?

Suppose that $G$ is finite.
Any configuration $\om\in\Om$ ($=\{ 0,1\}^E$) gives rise to a
dual configuration $\om^\tD $ lying in the space $\Om^\tD 
 =\{ 0,1\}^{E^\tD }$
defined by
$\om^\tD (e^\tD)=1-\om (e)$. As before, to each configuration 
$\om^\tD $ corresponds the
set $\eta (\om^\tD )=\{ e^\tD \in E^\tD : 
\om^\tD (e^\tD )=1\}$ of its `open
edges'. 
Let $f(\om )$ be the number of faces of the graph $(V,\eta (\om ))$,
including the infinite face
By drawing a picture, one may easily be convinced
(see Figure 6.1) that the faces of 
$(V,\eta (\om ))$ are in one--one correspondence with the components of
$(V^\tD ,\eta (\om^\tD ))$,
and therefore
$f(\om ) =k(\om^\tD )$,
in the obvious notation.
We shall make use of Euler's formula (see
[\ci{Wil}]),
$$
k(\om )= |V|-|\eta (\om )| +f(\om )-1, \qq\om\in\Omega.
$$
 
\topinsert
\centerline{\epsfxsize=7cm \epsfbox{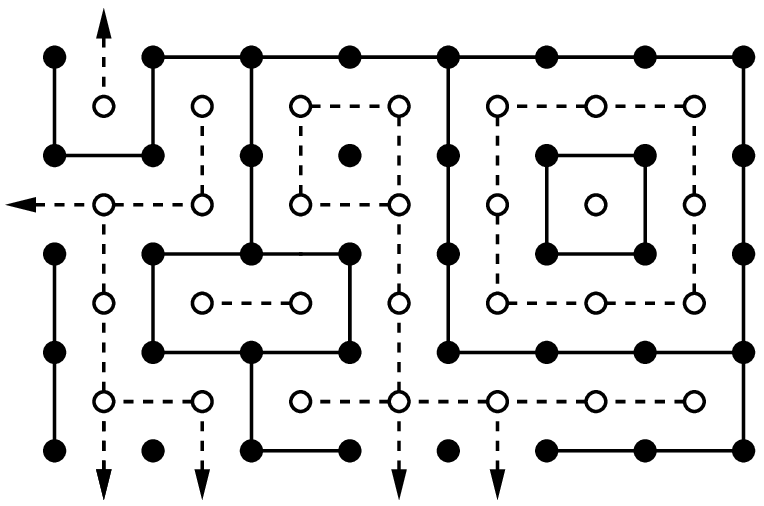}}
\capt{6.1}{A primal configuration $\om$ (with solid lines and vertices)
and its dual
configuration $\om^{\text{d}}$ (with dashed lines and hollow vertices).
The arrows join the given vertices of the dual
to a dual vertex in the infinite face.  Note that each face of the primal
graph (including the `infinite face')
corresponds to a unique component of the dual graph.}
\endinsert

The random-cluster measure on $G$ is given by
$$
\phi_{G,p,q}(\om )\propto\left(\frac{p}{1-p}\right)^{|\eta (\om )|}
q^{k(\om )}, \qq \om\in\Om.
$$
Using Euler's formula and the equality $f(\om)=k(\om^\tD)$, we find that
$$
\phi_{G,p,q}(\om )=\phi_{G^\td ,p^\td ,q}(\om^\td )\qq\text{for
$\om\in\Om$},
$$
where the dual parameter $p^\td$ is given according to
$$
\frac{p^\td}{1-p^\td}=\frac{q(1-p)}{p}.
$$
The unique fixed point of the mapping $p\mapsto p^\td$ is
given by $p=\psd q$ where
$\psd q={\sqrt q}/(1+\sqrt q)$.
We note at this point that
$$
\phi_{G,\psd q,q}(\om)\propto q^{\frac12|\eta(\om)|+k(\om)}
\propto q^{\frac12(k(\om^\tD)+k(\om))},
$$
by Euler's formula. This representation of the \rc\ measure at the `self-dual
point' $\psd q$ highlights the duality of measures.

Turning to the square lattice, let $\La_n=[0,n]^2$, whose dual graph
$\La_n^\td$ may be obtained from $[-1,n]^2+(\tfrac12 ,\tfrac12 )$ by
identifying all boundary vertices. By the above,
$$
\phi_{\La_n ,p,q}^0(\om )=\phi_{\La_n^\td ,p^\td ,q}^1 (\om^\td )
$$
for configurations $\om$ on $\La_n$ (and with a small `fix' on the
boundary of $\La_n^\td$). Letting $n\to\infty$, we obtain that
$\fpq^0(A)=\phi_{p^\td ,q}^1(A^\td )$
for all cylinder events $A$, where $A^\td =\{\om^\td :\om\in A\}$.

\subsection{6.2}{Value of the critical point}
Consider the random-cluster process on the two-dimensional lattice
$\LL^2=(\Z^2 ,\EE^2)$, with parameters $p$ and $q$ satisfying $q\geq
1$. The following remarkable conjecture is widely believed \OP.

\proclaim{Conjecture 6.\Thm} Let $q\ge 1$.
The critical value $\pc(q)$ of $\LL^2$ is
given by
$$
\pc(q) = \frac{\sqrt q}{1+\sqrt q}\qq\text{for } q\geq 1.
$$
\endproclaim
\label{exactpc}

This conjecture is known to hold when $q=1$, $q=2$,
and for $q\ge 25.72$.
The $q=1$ case was answered by Kesten [\ci{Kes}]
in his famous proof that the critical probability of bond percolation on
$\LL^2$ is $\tfrac12$. For $q=2$, the value of $\pc(2)$ given above
agrees with the celebrated calculation by Onsager [\ci{Ons}] of the
critical temperature of the Ising model on $\Z^2$,
and is implied by probabilistic results in the modern vernacular of
[\ci{ABF}].
The formula for $\pc(q)$ has been established rigorously in [\ci{LMR},
\ci{LMMRS}] for sufficiently large (real) values of $q$, specifically
$q\ge 25.72$
(see also [\ci{G96b}]).

Conjecture 6.\refer{exactpc} arises in a
natural manner from the observation that $\LL^2$
is a self-dual graph,
and $p=\psd q=\sqrt q/(1+\sqrt q)$ is the self-dual point of
a \rc\ measure on $\LL^2$ with parameters $p$, $q$.

Several other remarkable conjectures about the phase transition in $\LL^2$
may be found in the physics literature (see [\ci{Bax}]), as consequences
of `exact' but non-rigorous arguments involving ice-type models.
These include exact formulae for the asymptotic behaviour
of the partition function $\lim_{\La\uparrow\Z^2} \{Z_{\La,p,q}\}^{1/|\La|}$,
and also for the edge-densities $h^b(\psd q,q)
=\phi_{\psd q,q}^b(e\text{ is open})$, $b=0,1$, at the self-dual
point $\psd q$.

Progress towards a rigorous verification
of the conjecture may be summarised briefly as follows.
Using an argument, [\ci{Zha}],  taken from
percolation using the uniqueness of infinite
open clusters, we obtain by duality
that
$\t^0\bigl( \psd q,q\bigr) =0$
(see [\ci{G93}, \ci{Wel93}]), whence the critical value of the square lattice
satisfies
$\pc (q)\geq \psd q$ for $q\geq 1$. 
The complementary
inequality $\pc(q)\le \psd q$ has eluded mathematicians
despite progress by physicists, [\ci{HKW}].

\comment
The upper bound for $\pc(q)$ is obtainable via a crude and
non-rigorous argument. We may accept the following picture. If
$\psd q<\pc(q)$, then all components of the primal process 
with $p=\psd q$ are finite, and they
are islands which float in an infinite open ocean of the dual lattice.
Thus $p^\tD>\pc(q)$. However, $p^\tD=\psd q$, and thus we arrive
at a contradiction.  It follows that
$\pc(q)\le\psd q$ for $q\geq 1$. Note the fallacy in this
chain of deduction: there exist
configurations in which all primal and dual components are finite.
\endcomment

Suppose on the contrary that $\pc(q)>\psd q$, so that
$\pc(q)^\td<\psd q$. For $p\in (\pc(q)^\tD,\pc(q))$ 
we have also that $p^\tD\in (\pc(q)^\tD,\pc(q))$. Therefore, 
for $p\in(\pc(q)^\tD,\pc(q))$, both primal and dual 
processes comprise, almost surely,
the union of {\it finite\/} open clusters. This contradicts the intuitive
picture, supported
for $p\ne\pc(q)$ by our knowledge of percolation,
of finite clusters of one process floating in an infinite
ocean of the other process.

Exact values for the critical points of the triangular
and hexagonal lattices may be conjectured similarly, using graphical duality
together with the star--triangle transformation, [\ci{Bax}, \ci{KJ}].

Rigorous numerical upper
bounds of impressive accuracy
have been achieved for the square lattice and other two-dimensional
lattices via an extension of the basic model to a larger
class termed in [\ci{Al01}] the `asymmetric \rc\ model'.
The bound in question for $\LL^2$ is
$$
\pc(q)\le \frac{\sqrt q}{\sqrt{1-q^{-1}} +\sqrt q},\qq q\ge 1.
$$
For example, when $q=10$, we have that $0.760 \le \pc(10)\le 0.769$,
to be compared with the conjecture that $\pc(10)=\sqrt{10}/(1+\sqrt{10})
\simeq 0.760$. A valuable consequence of the comparison methods
developed in [\ci{Al01}] is the exponential decay of connectivity
functions when $q>2$ and $p$ is such that
$$
p<\psd {q-1}=\frac{\sqrt{q-1}}{1+\sqrt{q-1}}.
$$

\subsection{6.3}{First-order phase transition}
There is a special argument discovered first for Potts models, [\ci{LMR}],
which may be used to show first-order phase transition
when $q$ is sufficiently large.

Let $a_n$ be the number of self-avoiding walks on $\LL^2$ 
beginning at the origin. It is standard, [\ci{MS}], that $a_n^{1/n}\to\mu$ 
as $n\to\infty$,
for some constant $\mu$ called the {\it connective constant\/} of the lattice.
Let $Q= \bigl\{\frac12\bigl(\mu + \sqrt{\mu^2 - 4}\bigr)\bigr\}^4$.
We have that $2.620 < \mu < 2.696$ (see [\ci{Sl95}]),
whence $21.61 <  Q < 25.72$. We set
$$
\psi(q) = \frac 1{24} \log\left\{\frac{(1+\sqrt q)^4}{q\mu^4}\right\},
$$
noting that $\psi(q) > 0$ if and only if $q > Q$.
We write $B(n)=[-n,n]^2$.

\proclaim{Theorem 6.\Thm{} [\ci{G96b}, \ci{LMR}]} If $d=2$ and $q > Q$ then
the following hold.
\ii{\rm (a)} The critical point is given by $\pc (q)=\sqrt q/(1+\sqrt
q)$.
\ii{\rm (b)} We have that $\t^1(\pc (q),q)>0$.
\ii{\rm (c)} For any $\psi <\psi (q)$ and all large $n$,
$\phi_{\pc (q),q}^0 \bigl( 0\lra\pd B(n)\bigr)\leq e^{-n\psi}$.
Hence, in particular, $\t^0(\pc(q),q)=0$.
\endproclaim
\label{first2d}

The idea of the proof is as follows. 
There is 
a partial order on circuits of $\LL^2$ given
by $\Gamma\le \Gamma'$ if the bounded component of $\R^2\setminus\Gamma$
is a subset of the bounded component of $\R^2\setminus \Gamma'$.
We work at the self-dual
point $p=\psd q$, and with the box $B(n)$ with
wired boundary conditions. 
An `outer contour' is defined to be a circuit
$\Gamma$ of the dual graph $B(n)^\td$ all of whose edges are open
in the dual (that is, they traverse closed edges in the primal graph $B(n)$),
and which is maximal with this property.
Using self-duality, one may show that 
$$
\phi_{B(n),\psd q,q}^1(\Gamma \text{ is an outer circuit}) \le 
\frac 1q\left(\frac q{(1+\sqrt q)^4}\right)^{|\Gamma|/4},
$$
for any given circuit $\Gamma$ of $B(n)^\td$.
Combined with a circuit-counting argument of Peierls-type
involving the connective constant, this estimate implies  
after a little work the claims of
Theorem 6.\refer{first2d}. The idea of the proof 
appeared in [\ci{LMR}] in the context of
Potts models, and the \rc\ formulation may be found
in [\ci{G96b}].

We stress that corresponding 
conclusions may be obtained for general $d$ ($\geq 2$)
when $q$ is sufficiently large ($q > Q(d)$ for suitable $Q(d)$),
as shown in [\ci{LMMRS}] using so-called
Pirogov--Sinai theory. Whereas, in the case $d=2$,
the above duality provides an especially simple
proof, the proof for general $d$ utilises nested sequences of 
surfaces of $\R^d$ and requires a control of 
the effective boundary conditions within
the surfaces.  

\subsection{6.4}{\SLE\ limit when $q\le 4$}
Many exact calculations are `known' for critical processes
in two dimensions, but the physical arguments involved have
sometimes appeared in varying degrees 
magical or revelationary to mathematicians.
The new technology of stochastic L\"owner evolutions
(\SLE),
discovered by Schramm [\ci{Sch00}] and mentioned in Section 5.3,
threatens to provide a rigorous underpinning of many such
arguments in a manner most consonant with modern probability theory.
Roughly speaking, the theory of \SLE\ informs us of the correct weak limit
of a critical process in the limit of large spatial scales, and in addition
provides a mechanism for performing calculations for the limit process.

Let $\HH=(-\oo,\oo)\times(0,\oo)$ be the upper half-plane of $\R^2$,
with closure $\oHH$.  We view $\HH$ and $\oHH$ as subsets of the complex plane.
Consider the ordinary differential equation
$$
\frac d{dt}g_t(z)=\frac 2{g_t(z)-B_{\kappa t}},\qquad z\in\oHH\setminus\{0\},
$$
subject to the boundary condition $g_0(z)=z$,
where $t\in[0,\oo)$, $\kappa$ is a positive constant,
and $(B_t:t\ge 0)$ is a standard Brownian motion.
The solution exists when $g_t(z)$ is bounded away from $B_{\kappa t}$.
More specifically, for $z\in \oHH$,
let $\tau_z$ be the infimum of all times $\tau$
such that 0 is a limit point of
$g_s(z)-B_{\kappa s}$ in the limit as $s\uparrow \tau$.
We let 
$$
H_t=\{z\in\HH: \tau_z > t\}, \qq K_t=\{z\in\oHH: \tau_z\le t\},
$$
so that $H_t$ is open, and $K_t$ is compact.
It may now be seen that $g_t$ is a conformal homeomorphism
from $H_t$ to $\HH$.

We call $(g_t:t\ge 0)$ a 
{\it stochastic L\"owner evolution\/} (\SLE)
with parameter $\kappa$, written \SLE$_\kappa$, and we call
the $K_t$ the {\it hulls\/} of
the process. 
There is good reason to believe that the family $K=(K_t:t\ge 0)$ 
provides the
correct scaling limit of a variety of random spatial processes,
the value of $\kappa$ being chosen according to the process 
in question. General properties of \SLE$_\kappa$, 
viewed as a function of $\kappa$
have been studied in [\ci{RS01}, \ci{W02}], and a beautiful
theory has emerged. For example, the hulls $K$ form almost surely a simple path
if and only if $\kappa\le 4$.
If $\kappa> 8$, then \SLE$_\kappa$ generates
almost surely a space-filling curve.

Schramm [\ci{Sch00}] has identified the relevant value of $\kappa$
for several different processes, and has indicated that
percolation has scaling limit \SLE$_6$,
but full rigorous proofs are incomplete. In the case of percolation,
Smirnov [\ci{Smi0}, \ci{Smi}] has proved the very remarkable result that,
{\it for site percolation on the triangular lattice\/},
the scaling limit exists and is \SLE$_6$ (this last statement is illustrated
and partly explained in Figure 6.2), but the existence of the limit
is open for other lattices and for bond percolation.

\topinsert
\figure
\centerline{\epsfxsize=12cm \epsfbox{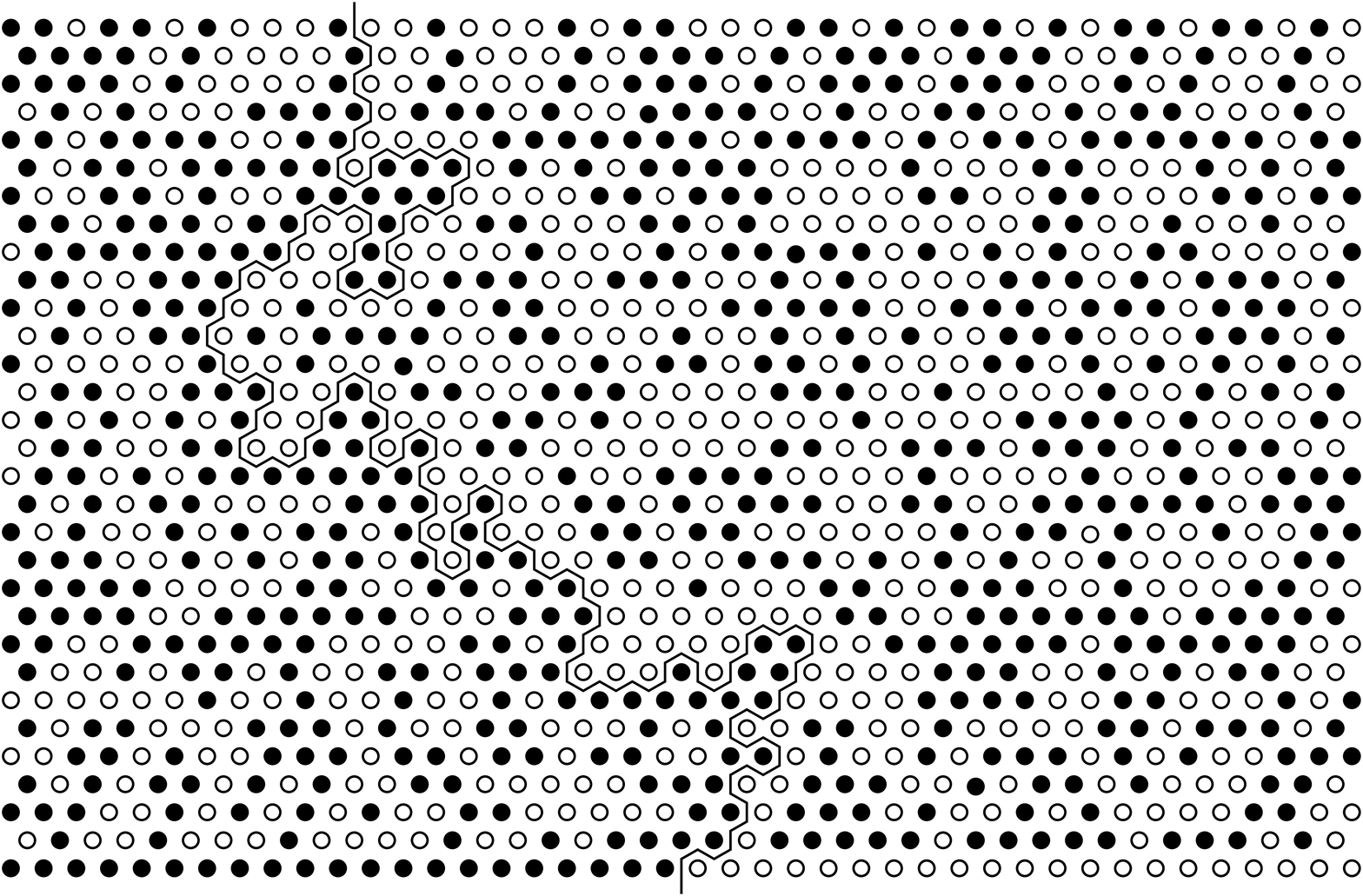}}

\capt{6.2}{Site percolation on the triangular lattice with $p$ equal
to the critical point $\frac12$, and with a mixed boundary condition
along the lower side. The interface
traces the boundary between the white and the black
clusters touching the boundary, and is termed the `exploration process'.
In the limit of small lattice-spacing,
the interface converges in
a certain manner to the graph of
a function which satisfies the L\"owner differential
equation driven by a Brownian motion with variance parameter
$\kappa=6$.}
\endinsert

It is possible to perform calculations on stochastic L\"owner evolutions,
and in particular to confirm, [\ci{LSW6}, \ci{SmiW}], the values of
many critical exponents associated with percolation
(or, at least, site percolation on the triangular lattice).
The consequences are in agreement with predictions of
mathematical physicists previously considered near-miraculous
(see [\ci{G99}], Chapter 9). In addition, \SLE$_6$ satisfies
the appropriate version of Cardy's formula, [\ci{Car}, \ci{LPS}].

The technology of \SLE\ is a major piece of contemporary
mathematics which promises to explain phase transitions in
an important class of two-dimensional disordered systems,
and to help bridge
the gap between probability theory and conformal field theory.
It has already provided complete explanations of conjectures,
by mathematicians and physicists, associated with
two-dimensional Brownian motions and specifically their intersection exponents
and fractionality of frontier, [\ci{LSW5}, \ci{LSW4}].

Extra work is needed in order to prove the validity of the limiting
operation for other percolation models and random processes. 
In another remarkable recent paper [\ci{LSW7}],
Lawler, Schramm, and Werner have verified the existence of
the scaling limit for loop-erased random
walk and for the uniform spanning tree Peano curve,
and have shown them to be \SLE$_2$ and \SLE$_8$ respectively.
It is believed that self-avoiding walk on $\LL^2$, [\ci{MS}], has
scaling limit \SLE$_{8/3}$.

We turn now to the \rc\ model on $\LL^2$ with
parameters $p$ and $q$. For $1\le q< 4$, it is believed that the percolation
probability $\theta(p,q)$, viewed as a function of $p$,
is continuous at the critical
point $\pc(q)$ \OP, and furthermore that $\pc(q)=\sqrt q/(1+\sqrt q)$.
It seems likely that, when re-scaled in the manner similar to
that of percolation (illustrated in Figure 6.2), the exploration
process of the model converges to a limit process
of \SLE\ type. It then remains only to specify the parameter $\kappa$
of the limit in terms of $q$. It has been
conjectured  in [\ci{RS01}] that 
$\kappa$ satisfies $\cos(4\pi/\kappa)=-\frac12\sqrt q$, $\kappa\in(4,8)$.
This value is consistent with Smirnov's theorem [\ci{Smi}],
and also
with the finding of [\ci{LSW7}] that
the scaling limit of the uniform
spanning tree Peano curve is \SLE$_8$, on recalling that the uniform
spanning tree measure is obtainable as a limit
of the \rc\ measure as $p,q\downarrow 0$ (recall Section 2.4).

There are uncertainties over how this programme will develop.
For a start, the theory of \rc\ models is not as complete
as that of percolation and of the uniform spanning tree. Secondly,
the existence of spatial limits is currently known only
in certain special cases. The programme is however ambitious and full
of promise, and should ultimately
yield a full picture of the critical
behaviour --- including values of exponents --- of \rc\ models,
and hence of Ising/Potts models, with $q\le 4$.

\endsection

\section{7. On complete graphs and trees}
While considerations of `real space--time' support the study
of such models on lattices such as $\LL^d$, 
it has proved rewarding also to analyse
the \rc\ model on certain other graphs including complete graphs and trees.
It is the presence of circuits in the underlying graph which is the root
cause of dependence between the states of
edges, and for this reason it is the complete
graph which provides an appropriate setting for what is 
termed `mean-field theory', in  which
vertices `interact' with {\it all\/} other vertices rather than with a selected subset
of `neighbours'. Trees, on the other hand, contain no circuits, and
their random-cluster theory is therefore sterile unless
one introduces boundary conditions. [A different approach to mean-field
theory has been studied in [\ci{KS}], namely
on $\LL^d$ for large $d$.]

\subsection{7.1}{On complete graphs}
The mean-field Potts model may be formulated as a Potts model on the
complete graph $K_n$, being the graph with $n$ labelled vertices every
pair of which is joined by an edge. The study of such a process dates
back at least to 1954, [\ci{KMS}],
and has been continued over the last fifty years
[\ci{BGJ}, \ci{KS}, \ci{Wu}]. The model is exactly soluble in the sense that
quantities of interest may be calculated exactly
and rigorously. It is therefore not
surprising that the corresponding random-cluster models (for real
$q$) have `exact solutions' also, [\ci{BGJ}]. 

Consider the random-cluster measure $\psi_{n,\lambda ,q}=
\phi_{K_n,\lambda/n,q}$ on the
complete graph $K_n$, having parameters $p=\lambda /n$ and $q$;
this is the appropriate scaling to allow an interesting limit
as $n\to\infty$. In the
case $q=1$, this measure is product measure, and therefore the ensuing
graph is an Erd\H os--R\'enyi random graph [\ci{Bol}, \ci{JLR}]. The overall
picture for general values of $q$ is rather richer than for the case
$q=1$, and many exact calculations may
be performed rigorously. It turns out that the phase transition is of first-order if and
only if $q>2$, and the behaviour of the system depends on how
$\lambda$ compares with a `critical value' $\lac(q)$ taking the
value
$$
\lac(q)=\cases
q&\text{if $0<q\leq 2$,}\\
2\left(\dfrac{q-1}{q-2}\right)\log (q-1)&\text{if $q>2.$}
\endcases
$$

\topinsert
\figure
\mletter{\theta(\lambda,q)}{-0.15}{.25}
\mletter{\theta(\lambda,q)}{4.05}{.25}
\mletter{\theta(\lambda,q)}{8.2}{.25}
\mletter{\lac(q)}{1.3}{2.9}
\mletter{\lac(q)}{5.5}{2.9}
\mletter{\lac(q)}{9.7}{2.9}
\mletter{\lambda}{3.6}{2.9}
\mletter{\lambda}{7.8}{2.9}
\lastletter{\lambda}{12}{2.9}
\centerline{\epsfxsize=12cm \epsfbox{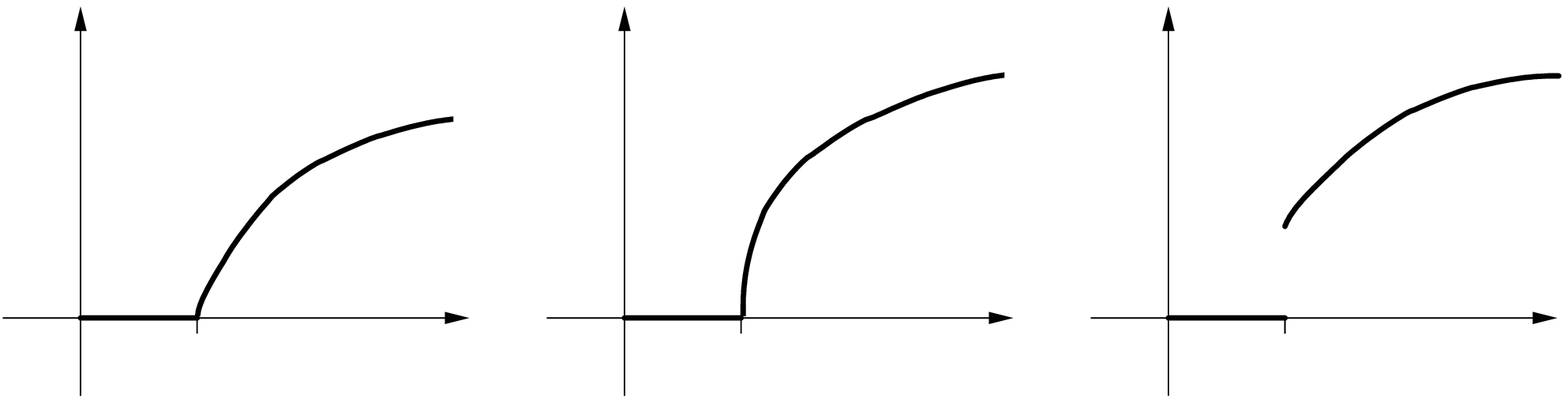}}
\line{\phantom{c}\hfil}
\line{\hfil $q<2$ \hskip 3.4cm $q=2$ \hskip 3.4cm $q>2$ \hfil}

\capt{7.1}{The function $\t(\lambda,q)$ for the three cases $q<2$, 
$q=2$, $q>2$.}
\endinsert

From the detailed picture described in [\ci{BGJ}]  the
following information may be extracted. The given properties occur with
$\psi_{n,\lambda ,q}$-probability tending to 1 as $n\to\infty$.
\ssubsection{I}{Subcritical case, when $\lambda <\lac(q)$.} 
The largest component of
the graph is of order $\log n$.
\ssubsection{II}{Supercritical case, when $\lambda >\lac(q)$.} 
There is a `giant
component' having order $\t (\lambda ,q)n$ where $\t$ is defined to be
the largest root of the equation
$$
e^{\lambda\t} = \frac{1+(q-1)\t}{1-\t}.
$$
\ssubsection{III}{Critical case, when $\lambda =\lac(q)$, $0<q\leq 2$.} 
The largest
component has order $n^{2/3}$.
\ssubsection{IV}{Critical case, when $\lambda =\lac(q)$, $q>2$.} The largest
component is either of order $\log n$ or of order $\t (\lambda ,q)n$,
where $\t$ is given as in case II above.

The dichotomy between first- and second-order phase transition is seen
by studying the function $\t (\lambda ,q)$, sketched in Figure 7.1. When
$0<q\leq 2$, the function $\t (\lambda ,q)$ descends continuously to 0 as
$\lambda\down\lac(q)$. On the other hand, this limit is strictly
positive when $q>2$.

The above results are obtained via a relationship between the model
for general $q$ and the model for the special case $q=1$.
The latter system 
has been analysed extensively, [\ci{Bol}, \ci{JLR}].
We illustrate the argument in the case $q \ge 1$; a similar approach is valid
when $q<1$. Consider the open clusters $C_1,C_2,\dots,C_m$ of 
a sample from the random-cluster
measure $\phi_{K_n,p,q}$. We colour each such cluster {\it red\/}
with probability $\rho$, and {\it white\/} otherwise, different
clusters receiving independent colours. We delete all vertices in white
clusters, and let $H$ denote the remaining graph, comprising a certain
random number $N$ of vertices (from the red clusters) together with
certain open edges joining pairs of  them. It may be seen that,
conditional on the value of $N$, the measure governing $H$ is the \rc\ measure
with parameters $p$ and $q\rho$. We choose $\rho=1/q$ to obtain an Erd\H os--R\'enyi
random graph on a {\it random\/} set of vertices. This is the
observation which permits the full analysis to proceed.

One consequence of this study is an explicit
identification of the exponential asymptotics
of the partition function $Z_{K_n,\lambda/n,q}$, namely of the limit function
$$
f(\lambda,q)=\lim_{n\to\infty}\left\{\frac1n\log Z_{K_n,\lambda/n,q}\right\}. 
$$
This provides
information via the G\"artner--Ellis theorem, [\ci{DZ}],
concerning the large-deviation theory
of the number of clusters in such systems.

\subsection{7.2}{On trees and non-amenable graphs}
Whereas physical considerations support the study of interacting systems
on finite-dimensional lattices, 
mathematicians have been drawn also to the study of general graphs,
thus enabling a clearer elucidation of the mathematical
structure of such systems (see the discussion in [\ci{Sch01}]).
A subject of special focus has been the class of graphs
for which the ratio of surface
to volume of finite boxes does not approach zero in the limit as the
size of the box tends to infinity. 
A prime example of such a graph is an infinite
regular tree with vertex degree at least three.
We make the distinction more concrete as follows.  Let $G=(V,E)$
be an infinite connected graph with
finite vertex degrees. For $W\subseteq V$, we define its {\it boundary\/}
$\partial W$ to be the set of all $w\in W$ having some neighbour $v$
not belonging to $W$. The 
countably infinite graph $G$ is called {\it amenable\/} if
its `Cheeger constant'
$$
\kappa(G)=\inf\left\{\dfrac{|\partial W|}{|W|}: W\subseteq V,\ 0<|W|<\infty\right\}
$$
satisfies $\kappa(G)=0$; $G$ is called {\it non-amenable\/} if $\kappa(G)>0$.
It is easily seen that $\LL^d$ is amenable, whereas an 
infinite regular tree with degree at least three is non-amenable.

The role of amenability in probability theory has been evident since
the work of Kesten [\ci{K59a}, \ci{K59b}] 
concerning random walks on a general graph $G$.
More relevant to this review are [\ci{BS}, \ci{GN}], 
which consider the number of infinite
clusters in the bond percolation model on $G$. Suppose $G$ is
a quasi-transitive graph (that is, its vertex set has
only finitely many orbits under its automorphism group).
Suppose
in addition that $G$ is amenable. Consider
bond percolation on $G$ with density $p$. 
It may be proved as in [\ci{BK}, \ci{GKN}] that
the number $I$ of infinite
open clusters satisfies
$$
\text{either}\q P_p(I=0)=1\q\text{or}\q P_p(I=1)=1.
$$
That is, if an infinite open cluster
exists, then it is almost surely unique. Under similar
assumptions on a non-amenable graph $G$,
it is believed but not yet proved in full generality
that there exists an interval of values of $p$ for which
$P_p(I=\infty)=1$; see, for example, the discussion in [\ci{L00}].
A corresponding question for \rc\ models is to ascertain for which
graphs $G$ and values of $q$ there is non-uniqueness of \rc\ measures
for an interval of values of $p$. 
[Recall Theorem 3.\refer{countable}, easily extended to
more general amenable graphs, which states that, for $q\ge 1$,
there is a unique \rc\ measure on
$\LL^d$ for all except at most countably many values of $p$.]
See [\ci{HJL1}, \ci{Jon99}] and especially [\ci{Sch01}, 
Section 6.1] for recent accounts
of this and associated questions, and [\ci{Hag96}] for
an analysis of \rc\ measures on regular trees.
\endsection

\section{8. Time-evolutions of \rc\ models}
Let $\mu$ be a probability measure on a space $(\Om,\sF)$. We may
study stochastic processes taking values in $\Om$ which converge weakly
to $\mu$ in the limit of large times.
There are a multiplicity of reasons for and benefits in studying
time-evolutions. First, physical systems generally have dynamics as well
as equilibria. Secondly, new questions of interest arise, such as
that of the estimation of a relaxation time. Thirdly, 
the dynamics
thus introduced can yield a new technique for studying
the limit measure $\mu$. 

When studying a physical
system, it is often acceptable to restrict oneself to
dynamics which are reversible in time.
In Section 8.1, we  describe a natural reversible dynamic
for a \rc\ model, akin to the Glauber dynamics of the Ising model.
This dynamic permits an extension which couples together
the \rc\ measures on a given graph as $p$ and $q$ range
over their possible values.

The problem commonly arises in statistics, computer science, and
statistical physics of how to obtain a sample from a system
governed by a probability measure with complex structure.
In Section 8.2 we summarise the Propp--Wilson
`coupling from the past'
approach, [\ci{PW96}], to this problem in the context of the \rc\ measure.

Since Potts models may be obtained from \rc\ models, there is
an interplay between the dynamics for these two systems. A famous 
instance of this relationship is the so-called Swendsen--Wang dynamic
[\ci{SW}], which is described in Section 8.3.

We assume throughout this section that $G=(V,E)$ is a finite
connected graph, and that $\fpq$ is the \rc\ measure on $\Om=\{0,1\}^E$
with $0<p<1$ and $q>0$.

\subsection{8.1}{Reversible dynamics}
We shall consider transitions from a configuration $\om$ to configurations
which differ from $\om$ on one edge only. Thus we introduce the following
notation. For $\om\in\Om$ and $e\in E$, let $\om^e$ and $\om_e$ be the
configurations obtained by `switching $e$ on' and `switching $e$
off', respectively, that is
$$
\om^e(f)=\cases 1 &\text{if } f=e,\\ \om(f) &\text{if } f\ne e,\endcases
\qquad
\om_e(f)=\cases 0 &\text{if } f=e,\\ \om(f) &\text{if } f\ne e.\endcases
$$

Let $(X_t: t\ge 0)$ be a Markov chain, [{\ci{GS}],
on the state space $\Om$ with generator $Q=\{q_{\om,\om'}: \om,\om'\in\Om\}$
satisfying
$$
q_{\om_e,\om^e} = p, \quad q_{\om^e,\om_e}=(1-p)q^{D(e,\om_e)},
\qquad\text{for } \om\in\Om,\ e\in E,
$$
where $D(e,\xi)$ is the indicator function of the event that the endpoints of
$e$ are joined by no open path of $\xi$. This specifies the rate at which 
single edges
are acquired or lost. We set $q_{\om,\xi}=0$ if $\om$ and $\xi$
differ on two or more edges, and we choose the diagonal elements $q_{\om,\om}$
in such a way that $Q$, when viewed as a matrix, has row sums zero,
that is,
$$
q_{\om,\om}=-\sum_{\xi:\xi\ne\om} q_{\om,\xi}.
$$

It is elementary that the `detailed balance equations'
$$
\fpq(\om)q_{\om,\om'}=\fpq(\om')q_{\om',\om}, \q \om,\om'\in \Om,
$$
hold,
whence $X$ is reversible with respect to
$\fpq$. It follows by the irreducibility of the chain that
$X_t\Rightarrow \fpq$ as $t\to\infty$ (where `$\Rightarrow$' denotes
weak convergence). There are of course many Markov
chains with generators satisfying the above detailed balance
equations, the important quantity is the ratio $q_{\om,\om'}/q_{\om',\om}$.

Two extensions of this dynamical structure which have proved
useful are as follows. The evolution may be specified in terms of
a `graphical representation' constructed via a family of
independent Poisson processes. This allows a natural coupling of
the measures $\fpq$ for different $p$ and $q$. Such couplings
are monotone in $p$ when $q\ge 1$. One may similarly couple
the unconditional measure 
$\fpq(\cdot)$ and the conditioned measure $\fpq(\cdot \mid A)$. 
Such couplings permit probabilistic
interpretations of differences of the form
$\phi_{p',q}(B\mid A) - \fpq(B)$ when $q\ge 1$, $p\le p'$,
and $A$ and $B$ are increasing, and this can be useful
in particular calculations (see [\ci{BGK}, \ci{G95a}, \ci{G93}]).

We turn now to the thermodynamic limit, and the question
of the structure of a Markovian \rc\ process on an infinite connected
graph. In the case $q\ge 1$, the above couplings are monotone 
in the choice of the underlying graph $G$. Therefore there 
exist `limit dynamics'
as $G$ passes through an increasing sequence of finite graphs.
Boundary conditions may be introduced, and one may obtain thereby
a certain Markov process $\zeta=(\zeta_t: t\ge 0)$ 
on the state space $[0,1]^{\EE}$,
where $\EE$ is the limiting (infinite) edge set. This process,
which does not generally have the Feller property,
generates a pair of `level-set processes' taking values
in $\{0,1\}^\EE$, defined for $0\le p\le 1$ by
$$
\zeta_t^{p,-}(e) = 1_{\{\zeta_t(e)>1-p\}},\q
\zeta_t^{p,+}(e) = 1_{\{\zeta_t(e)\ge 1-p\}},
\qq e\in \EE,
$$
where, as before, $1_A$ denotes the indicator function of
an event $A$. These two processes are Markovian
and are reversible with respect to the 
infinite-volume free and
wired random-cluster measures, respectively. See [\ci{G93}].

Note that the generator of the Markov chain given above
depends on the random variable $D(e,\om_e)$, and that this random variable
is a `non-local'
function of the configuration $\om$ in the sense that there is no absolute
bound on the distance from $e$ of edges whose states may be relevant 
to its value. It is this feature of non-locality which leads
to interesting complications linked in part to the
$0/1$-infinite-cluster property introduced before
Theorem 3.\refer{rcexist}. Further discussion
may be found in [\ci{G93}, \ci{PV}].

\subsection{8.2}{Coupling from the past}
In running Monte Carlo experiments, one requires the 
ability to sample from the probability measure $\fpq$. The Markov chain $X_t$
of Section 8.1 certainly converges weakly to $\fpq$ as $t\to\infty$, but
this is not as good as having in the hand
a sample with the {\it exact\/} distribution.
Random-cluster measures are well suited to the Propp--Wilson
approach to sampling 
termed `coupling from the past', [\ci{PW96}], and we sketch this
here. Some illustrations may be found in [\ci{Je98}].

First we provide ourselves with a discrete-time reversible Markov chain
$(Z_n: n\ge 0)$ on the state space
$\Om$ having invariant measure $\fpq$. 
The so-called {\it heat-bath algorithm\/} provides a suitable example 
of such a chain, and
proceeds as follows. At each stage, we pick a random edge $e$, chosen
uniformly from $E$ and independently of all earlier choices,
and we make $e$ open with the correct conditional probability,
given the configuration on the other edges.
The corresponding transition matrix
is given by $\Pi=\{\pi_{\om,\om'}: \om,\om'\in\Om\}$ where
$$
\align
\pi_{\om_e,\om^e} &=  \frac1{|E|}\cdot\frac {\fpq(\om^e)}{\fpq(\om^e)
  +\fpq(\om_e)},\\
\pi_{\om^e,\om_e} &=\frac1{|E|}\cdot\frac {\fpq(\om_e)}{\fpq(\om^e)
    +\fpq(\om_e)}.
\endalign
$$
A neat way to do this is as follows. Suppose that $Z_n=\om$.
Let $e_n$ be a random edge of
$E$, and let $U_n$ be uniformly distributed on the interval $[0,1]$, 
these variables being chosen independently of all earlier
choices. We obtain $Z_{n+1}$ from $\om$ by retaining the
states of all edges except possibly that of $e_n$.
We set 
$$
Z_{n+1}(e_n)=0\q \text{if and only if}\q
U_n\le  \frac {\fpq(\om_e)}{\fpq(\om^e)
    +\fpq(\om_e)}.
$$
Thus the evolution of the chain is determined
by the sequences $e_n$, $U_n$, and the initial state $Z_0$.
One may make this construction explicit by writing
$Z_{n+1}=\psi(Z_n, e_n, U_n)$ ($=\psi(\om, e_n, U_n)$) for some 
function $\psi:\Om\times E\times [0,1]\to\Om$.
It is easily seen by the Holley
condition of Section 3.1 that, if $q\ge 1$, and for every $e$ and $u$,
the function $\psi(\cdot,e,u)$ is non-decreasing
in its first argument.
It follows that the coupling is `monotone' in the sense that, if
$\om\le\om'$, then the chain starting at $\om$ lies at all times
beneath the chain starting at $\om'$ (using the partial order on $\Om$).

We let $W=(W(\om):\om\in\Om)$ be a vector
of random variables such that $W(\om)$ has the distribution
of $Z_1$ conditional on $Z_0=\om$. Following the scheme
described above, we may take $W(\om)=\psi(\om,e,U)$ where $e$ and $U$
are chosen at random.
Let $W_{-m}$, $m\ge 1$, be 
independent random vectors distributed as $W$, that is,
$W_{-m}(\cdot)=\psi(\cdot,e_m,U_m)$ where 
the set $\{(e_m,U_m):m\ge 1\}$ comprises independent pairs
of independent random variables, each $e_i$ being uniform on $E$,
and each $U_i$ being uniform on $[0,1]$.
We now construct a sequence $Y_{-n}$, $n\ge 1$, of random maps from 
$\Om$ to $\Om$ by the following inductive procedure.  First, for $\om\in\Om$,
we set $Y_{-1}(\om)=W_{-1}(\om)$.
Having found $Y_{-1}, Y_{-2},\dots,Y_{-m}$,
we define $Y_{-m-1}(\om)=Y_{-m}(W_{-m-1}(\om))$. That is, $Y_{-m-1}(\om)$
is obtained from $\om$ by passing in one step to $W_{-m-1}(\omega)$,
and then applying $Y_{-m}$ to this new state.
The exact dependence structure of
this scheme is an important ingredient of what follows.

We stop this process at the earliest time $m$ 
at which `coalescence has occurred',
that is, at the moment $M$ given by 
$M=\min\{m:Y_{-m}(\cdot) \text{ is the constant function}\}$.
It is a theorem, [\ci{PW96}], that $M$ is $\fpq$-a.s.\ 
finite and, for any $\om$,
the random output $Y_{-M}(\om)$ 
is governed exactly by the probability measure $\fpq$.

This procedure looks unwieldy, since $\Om$ is typically rather large,
but the reality is simpler when $q\ge 1$.
By the monotonicity of the above coupling when $q\ge 1$, 
it suffices
to follow the trajectories of the `smallest' and `largest' configurations,
namely those beginning, respectively,
with every edge closed and with every edge open.
The processes
starting at intermediate configurations remain sandwiched
between the extremal processes, for all times $t$. Thus one may
define $M$ by $M=\min\{m: Y_{-m}(0) = Y_{-m}(1)\}$,
where $0$ and $1$ denote the vectors of zeros and ones as before.

\comment
Samples observed by this method are provided in Figure 2.1.
\endcomment

\subsection{8.3}{Swendsen--Wang dynamics}
It is a major target of statistical physics to understand the time-evolution
of disordered systems, and a prime example lies in the
study of the Ising model.
A multiplicity of types of dynamics have been 
proposed. The majority of these
share a quality of `locality' in the sense that the evolution involves
changes to the states of vertices in close proximity to one another,
perhaps single spin-flips, or spin-exchanges.  The state space
is generally large, of size $2^N$ where $N$ is the number of vertices,
and the Hamiltonian has complicated structure. When subjected 
to `local dynamics', the process may approach equilibrium very slowly
(see [\ci{Mart97}, \ci{Sch98}] for accounts of recent work of relevance).
`Non-local dynamics', on the other hand,  have the
potential to approach equilibrium faster, since they permit
large jumps around the state space, relatively unconstrained by 
neighbourly relations. The \rc\ model has played a role in
the development of a simple but attractive such system,
namely that proposed by Swendsen and Wang [\ci{SW}] and described
as follows for the Potts model with $q$ states.

As usual, $G=(V,E)$ is a finite graph, typically a large box in $\Z^d$,
and $\Sigma = \{1,2,\dots,q\}^V$ is the state space of a Potts model on $G$.
We write $\Omega=\{0,1\}^E$.
Suppose that, at some time $n$, we have obtained a configuration 
$\s_n$ ($\in\Sigma$).
We construct $\s_{n+1}$ as follows.  Let $p=1-e^{-\beta J}$
where $0<\beta J <\infty$.

\ii{I.} We let $\omega_n\in\Omega$ be given 
as follows. For $e=\la x,y\ra \in E$,
$$
\align
\text{if $\s_n(x)\ne \s_n(y)$, let } &\om_n(e)=0,\\
\text{if $\s_n(x)=\s_n(y)$, let } &\om_n(e)=\cases 1 &\text{with probability $p$},\\
   0 &\text{otherwise}, \endcases
\endalign
$$
different edges receiving independent states. The edge configuration $\om_n$
is carried forward to the next stage.

\ii{II.} To each cluster $C$ of the graph $(V,\eta(\om_n))$ we assign 
an integer chosen uniformly at random from the set $\{1,2,\dots,q\}$,
different clusters receiving independent labels. We let $\s_{n+1}(x)$
be the value thus assigned to the cluster containing the vertex $x$.

It may be checked that the Markov chain $(\s_n: n\ge 0)$ has as unique
invariant measure the Potts measure on $\Sigma$ with parameters
$\beta$ and $J$. (Recall paragraph (c) of Section 2.3.)

The Swendsen--Wang algorithm leads to samples which
generally converge to equilibrium faster than those 
defined via local dynamics. This is especially
evident in the `high $\beta$' (or `low temperature')
phase, for the following reason. Consider for example
the simulation of an Ising model on a finite 
box with free boundary conditions, and suppose
that the initial state is $+1$ at all vertices.
If $\beta$ is large, then local dynamics result in samples which 
remain close
to the `$+$ phase' for a very long time. Only after
a long wait will the process achieve an average magnetisation
close to 0. Swendsen--Wang dynamics, on the other hand, can achieve
large jumps in average magnetisation even in a single step, since
the spin allocated to a given large cluster of the corresponding \rc\ 
model is equally likely to be either of the two possibilities.
A rigorous analysis of rates of convergence
is however incomplete. It turns out that, {\it at\/} the
critical point, Swendsen--Wang dynamics approach equilibrium
only slowly, [\ci{BCFKTVV}].
A further discussion is available in [\ci{GHM}].

Algorithms of Swendsen--Wang type have been described for
other statistical mechanical models having graphical representations
of \rc{}-type; see [\ci{ChMa1}, \ci{ChMa2}]. Related work may be found
in [\ci{Wol89}].

\endsection

\bigskip\goodbreak\flushpar
{\bf Acknowledgements.} GRG recalls John Hammersley passing to him
in 1971 a copy of Fortuin's thesis [\ci{F71}] in which much of
the basic theory is developed.
Piet Kasteleyn kindly filled out the origins
of \rc\ models in two letters addressed to
GRG in November 1992.
The author acknowledges the opportunity
given by the Landau Center of the Hebrew University, Jerusalem,
to deliver a course of lectures on the \rc\ model
during July 2001. Harry Kesten kindly criticised a draft
of the work. Thanks are due to Malwina \L uczak 
for her contributions
to discussions on certain topics in this paper,
and to \'Agoston Pisztora for reading and commenting on parts of it.
The further suggestions of  Christian Borgs, 
Olle H\"aggstr\"om, Russell Lyons, Roberto
Schonmann, Oded Schramm, and Alan Sokal have been appreciated. The paper
was completed during a programme at the Isaac Newton Institute (Cambridge).
\endsection

\def\JSP{Journal of Statistical Phy\-sics\/}

\def\CMP{Communications in Math\-emat\-ical Ph\-ys\-ics\/}
\def\RMP{Reviews in Mathematical Physics\/}
\def\TPR{The Physical Review\/}
\def\PTRF{Probability Theory and Related Fields\/}
\def\AP{Annals of Probability\/}

\def\JMP{Journal of Math\-ematical Phys\-ics\/}
\def\JAP{Journal of Applied Probability\/}

\def\JPhysA{Journal of Phys\-ics A: Math\-ematical and General\/}
\def\JPhysC{Journal of Phys\-ics C: Solid State Physics\/}
\def\JPhysF{Journal of Phys\-ics F: Metal Physics\/}

\def\PCPS{Proceedings of the Cambridge Philosophical Society\/}

\def\CPC{Combinatorics, Probability, Computing\/}

\def\ZfW{Zeit\-schrift f\"ur Wahr\-schein\-lich\-keits\-theorie 
und Verwandte Gebiete\/}
\def\PRL{Physical Review Letters\/}
\def\RSA{Random Structures and Algorithms\/}
\def\SPA{Stochastic Processes and their Applications\/}
\def\AMS{American Mathematical Society}
\def\AIHP{Annales de l'Institut Henri Poincar\'e, Probabilit\'es et
Statistiques\/}
\def\TAMS{Transactions of the \AMS\/}
\def\MPRF{Markov Processes and Related Fields\/}

\def\PRE{Physical Review E\/}
\def\PRB{Physical Review B\/}

\def\ECP{Electronic Communications in Probability\/}
\def\MS{Mathematica Scandinavica\/}
\def\AMM{American Mathematical Monthly\/}
\def\EJP{Electronic Journal of Probability\/}
\def\IJM{Israel Journal of Mathematics\/}
\def\MRL{Mathematics Research Letters\/}
\def\JCP{Journal of Chemical Physics\/}
\def\CRAS{Comptes Rendus des S\'eances de l'Acad\'emie des Sciences. S\'erie I.
Math\'ematique\/}
\hyphenation{Sprin-ger}
\def\nref#1\endref{}

\enddocument\vfill\eject
\end